\documentclass[11pt]{article}

\usepackage{amsmath,amssymb,amsthm,amsfonts,color}
\usepackage{indentfirst}
\usepackage{graphicx}
\usepackage{enumitem}
\usepackage{verbatim}
\usepackage{mathrsfs}
\usepackage{tikz}
\usepackage[colorlinks,allcolors=red]{hyperref}
\usepackage{graphicx}

\newtheorem{theorem}{Theorem}[section]
\newtheorem{lemma}[theorem]{Lemma}

\newtheorem{corollary}[theorem]{Corollary}
\newtheorem{definition}[theorem]{Definition}
\newtheorem{remark}[theorem]{Remark}

\def\sqr#1#2{\vbox{\hrule height .#2pt
\hbox{\vrule width .#2pt height #1pt \kern #1pt
\vrule width .#2pt}\hrule height .#2pt }}
\def\square{\sqr74}
\def\endproof{\hphantom{MM}\hfill\llap{$\square$}\goodbreak}

\def\implies{\Longrightarrow}

\def\Hat{\widehat}

\def\begi{\begin{itemize}}
\def\endi{\end{itemize}}
\def\bega{\begin{array}}
\def\enda{\end{array}}

\def\forall{\hbox{for every~ }}
\def\ov{\overline}
\def\ds{\displaystyle}
\def\R{\mathbb{R}}

\def\L{{\bf L}}

\def\T{{\mathcal T}}
\def\B{{\mathcal B}}
\def\S{{\mathcal S}}

\def\A{{\mathcal A}}
\def\forall{\hbox{for all}~}

\def\vp{\varphi}

\def\diams{\diamondsuit}

\def\O{{\mathcal O}}

\def\U{{\mathcal U}}

\def\F{{\mathcal F}}
\def\C{{\mathcal C}}
\def\Tilde{\widetilde}

\def\bfe{{\bf e}}

\def\meas{\hbox{meas}}
\def\v{\vskip 1em}
\def\vs{\vskip 2em}

\def\ve{\varepsilon}

\def\bel{\begin{equation}\label}
\def\eeq{\end{equation}}

\begin{document}

\title{\bf On the Regularity of Optimal Dynamic Blocking Strategies}
\vs

\author{Alberto Bressan and  Maria Teresa Chiri\\
\,
\\
Department of Mathematics, Penn State University \\
University Park, Pa.~16802, USA.\\
\,
\\
e-mails: axb62@psu.edu, mxc6028@psu.edu.
}
\maketitle
\begin{abstract} 
The paper studies a dynamic blocking problem, motivated by a model of
optimal fire confinement.   While the fire can expand with unit speed in all directions, barriers are constructed in real time.
An optimal strategy is sought, minimizing the total value of the burned 
region, plus a construction cost.    It is well known that optimal barriers exists.
In general, they are a countable union of compact, connected, rectifiable sets.
The main result of the present paper shows that optimal barriers are nowhere dense. The proof relies on new estimates on the reachable sets and on 
optimal trajectories for the fire, solving 
a minimum time problem in the presence of obstacles.
\end{abstract}

 
\vs
\section{Introduction}
\label{s:1}
\setcounter{equation}{0}
We consider the dynamic blocking problem introduced in \cite{B}, for a model 
of wildfire propagation \cite{S}.
To restrict the spreading of the fire, it is assumed 
that a barrier can be constructed, in real time.
This could be a thin strip of land which is either 
soaked with water
poured down from a helicopter, 
or cleared from all
vegetation using a bulldozer, or sprayed with fire extinguisher 
by a team of firemen.  In all cases,
the  fire will not cross that particular strip of land. 
Here the key point is that the barrier is being constructed at the same time 
as the fire front is advancing.
 
In this setting,
a natural problem is to find 
the  best possible  strategy.   In other words,
we seek the
optimal location of the barriers, in order
to minimize:
\bel{1}\hbox{[total value of the burned area] + 
[total cost for constructing the barriers]}\eeq
among all barriers that can be constructed in real time.

We consider here the simplest situation where the fire initially burns on an open set  $R_0$, and 
propagates with unit speed in all directions.
We assume
\begi
\item[{\bf (A1)}]
{\it The initial set
$R_0\subset\R^2$ is open, bounded, nonempty, connected, with Lipschitz boundary
$\partial R_0$.}
\endi

If barriers are not present, for each $t\geq 0$ the
set $R(t)$ reached by the fire is defined as 
\bel{1.2}\bega{rl}R(t)&\doteq~
\bigl\{ x(t)\,;~~x(\cdot)~~\hbox{is 1-Lipschitz}\,,~~
x(0)\in R_0\bigr\}\\[3mm]
&=~\bigl\{ x\in \R^2\,;~~d(x,R_0)<~t\bigr\}.\enda\eeq
Here and in the sequel, by 1-Lipschitz we mean a function with Lipschitz constant 1.
Moreover,
 $d(x,\Omega)$ denotes the distance of a point $x$ to the set $\Omega\subset\R^2$, while $\langle\cdot,\cdot\rangle$ 
is the Euclidean inner product in $\R^2$.  The closure and the
boundary of $\Omega$ are denoted by
$\ov \Omega$ and $\partial\Omega$ respectively.
By  $B(x,r)$ we denote the open ball centered at
$x$ with radius $r$.  More generally, for $\Omega\subset\R^2$, $B(\Omega,r)
= \{x\,;~d(x,\Omega)<r\}$
denotes the open neighborhood of radius $r$ around $\Omega$.
Finally, $m_1, m_2$ denote the 1-dimensional and 2-dimensional 
Hausdorff measure, respectively.

Next, we assume that 
the spreading of the fire can be controlled by constructing a
barrier. 

\begin{definition}\label{d:11} A {\bf barrier} $\Gamma\subset\R^2$ is a disjoint union of countably many 
compact connected, rectifiable sets, with finite total length.
\end{definition}
Throughout the following, we write
\bel{Gi}
\Gamma~=~\bigcup_{i\geq 1} \Gamma_i\eeq
to denote a barrier, as a union of its compact, rectifiable, connected  components.

Intuitively, we think of 
a barrier as 
a family of curves in the plane, 
which the fire cannot cross.
When a barrier $\Gamma$ is in place, the  set
reached by the fire is reduced. This leads 
to the definition of the new reachable set
\bel{1.6}R^\Gamma(t)\doteq~\Big\{
x(t)\,;~~x(\cdot)~\hbox{is 1-Lipschitz}\,,~~
x(0)\in R_0\,,
\quad  x(\tau)\notin \Gamma~~ \hbox{for
all}~\tau\in [0,t]\Big\}\,.\eeq
Clearly, in this case the burned set will be somewhat smaller:  
$R^\Gamma(t)\subseteq R(t)$
for every $t\geq 0$.
Since in our model the barrier is constructed at the same time as the fire propagates, a restriction
on its length must  be imposed.
\begin{definition}\label{d:12}  
Given a construction speed $\sigma>1$,  we say that 
the barrier $\Gamma$ is {\bf admissible} if
\bel{ad1}
m_1\bigl(\Gamma\cap \ov{R^\Gamma(t)}\bigr)~\leq~\sigma t\qquad\qquad\forall t\geq 0.\eeq
\end{definition}

\begin{remark} {\rm 
 For each $t\geq 0$, the set 
$$\gamma(t)~\doteq~\Gamma\cap \ov{R^\Gamma(t)}$$ 
appearing in (\ref{ad1})
is the part of the barrier $\Gamma$ touched by the fire at time $t$.
This is the portion that actually needs to be put in place within
time $t$, in order to restrain the fire. The remaining portion
$\Gamma\setminus \gamma(t)$ can be constructed at a later time.
This motivates the above definition.  The equivalence between different
formulations of the dynamic blocking problem was proved in \cite{BW1}.
}\end{remark}

\v

Fire propagation can equivalently be described in terms
of the minimum time function 
\bel{TG}T^\Gamma(x)~\doteq~\inf\bigl\{t\geq 0\,;~~x\in \ov{R^\Gamma(t)}
\bigr\}\,.\eeq
{}From the definition, it follows that $T^\Gamma$ is lower semicontinuous.
We think of  $T^\Gamma(x)$ as 
the minimum time needed for the fire to reach the point $x$,
starting  from $R_0$ and without crossing the barrier.
Notice that  $T^\Gamma(x)=+\infty$ if the fire never reaches a neighborhood of
$x$.
In general, the minimal time function $T^\Gamma$ can be 
computed by solving  a Hamilton-Jacobi equation with obstacles, namely
\bel{4}|\nabla T(x)|~=~1\qquad x\notin \Gamma\,,\eeq
\bel{44} T(x) = 0
\quad\hbox{if }~ x\in R_0\,.\eeq
For a precise definition and properties of this solution, see \cite{DLR}.
We recall that $T^\Gamma$ is locally an SBV function \cite{AFP}.   
The set where it
has jumps is contained inside $\Gamma$.
If the function $T^\Gamma$ is known, we can then recover
the region $R^\Gamma(t)$ burned within time $t$ as
$$\ov {R^\Gamma(t)}~=~\bigl\{ x\in\R^2\,;~T^\Gamma(x)\leq t\bigr\}.$$

Two mathematical problems can now be formulated.

{\bf (BP)  Blocking Problem.}  
Given a bounded open set $R_0$, decide whether there exists an admissible barrier
$\Gamma$ 
such that the entire region burned by the fire  
\bel{RGI}R^\Gamma_\infty~\doteq~\bigcup_{t>0} R^\Gamma(t)\eeq
  is bounded.

\v

{\bf(OP) Optimization Problem.} Given an initial 
set $R_0$ and a constant $c_0\geq 0$, find an admissible barrier
$\Gamma$ which minimizes the total cost
\bel{cost}{ \cal J}(\Gamma)~\doteq~m_2\bigl(R^\Gamma_\infty\bigr) +c_0\, m_1(\Gamma).\eeq
 \v
\begin{remark} {\rm For a given initial domain $R_0$,
the set $R^\Gamma_\infty$ in (\ref{RGI})
burned by the fire
can be characterized as the union of all connected components
of $\R^2\setminus \Gamma$ which intersect $R_0$.
For any bounded open set $R_0$, it is known \cite{B, Breview,
BBFJ, BW2} that a blocking strategy exists
if the construction speed is $\sigma>2$, while it does not exist 
if $\sigma\leq 1$.   The existence of a blocking strategy 
for $\sigma\in \,]1,2]$ is a challenging open problem. 
See the 
review \cite{Breview} for a more comprehensive discussion.
}\end{remark}

In a very general setting,  the existence of an optimal 
barrier was proved in \cite{BDL, DLR}.
Under the  assumption that this optimal barrier is the union of finitely many 
Lipschitz arcs, various necessary conditions  were derived in \cite{B, BW3, W}.
Indeed, assuming Lipschitz regularity,  one can reformulate the problem in the classical 
setting of  the Calculus of Variations, or within the 
theory of optimal control \cite{BPi, 
Cesari, FR}. Necessary conditions for optimality
are thus obtained in terms of the  
Euler-Lagrange equations, or  by applying the Pontryagin Maximum Principle. 
For example, when the initial set $R_0$ is a circle and the construction speed is
$\sigma>2$, 
among all simple closed curves, it is known that 
the admissible barrier that encloses the smallest burned area
is the union of an arc of circumference and two logarithmic spirals~\cite{BW4}.
 
Unfortunately, 
the results in \cite{BDL, DLR} only provide the existence of an optimal barrier 
$\Gamma^*$ with 
the minimal regularity properties stated in Definition~\ref{d:11}.  Namely, 
we only know that 
$\Gamma^*$ is the  union of countably many 
compact, connected, rectifiable sets.
It remains an outstanding open problem to close this gap, establishing further
regularity properties of the optimal barrier, so that necessary conditions for optimality can then be applied.  
In the present paper we take a step in this direction.
Our main goal is to prove

\begin{theorem}\label{t:1} 
For the optimization problem {\bf (OP)}, any optimal barrier  
 $\Gamma$ is nowhere dense.
\end{theorem}

\begin{figure}[htbp]
   \centering
 \includegraphics[width=0.9\textwidth]{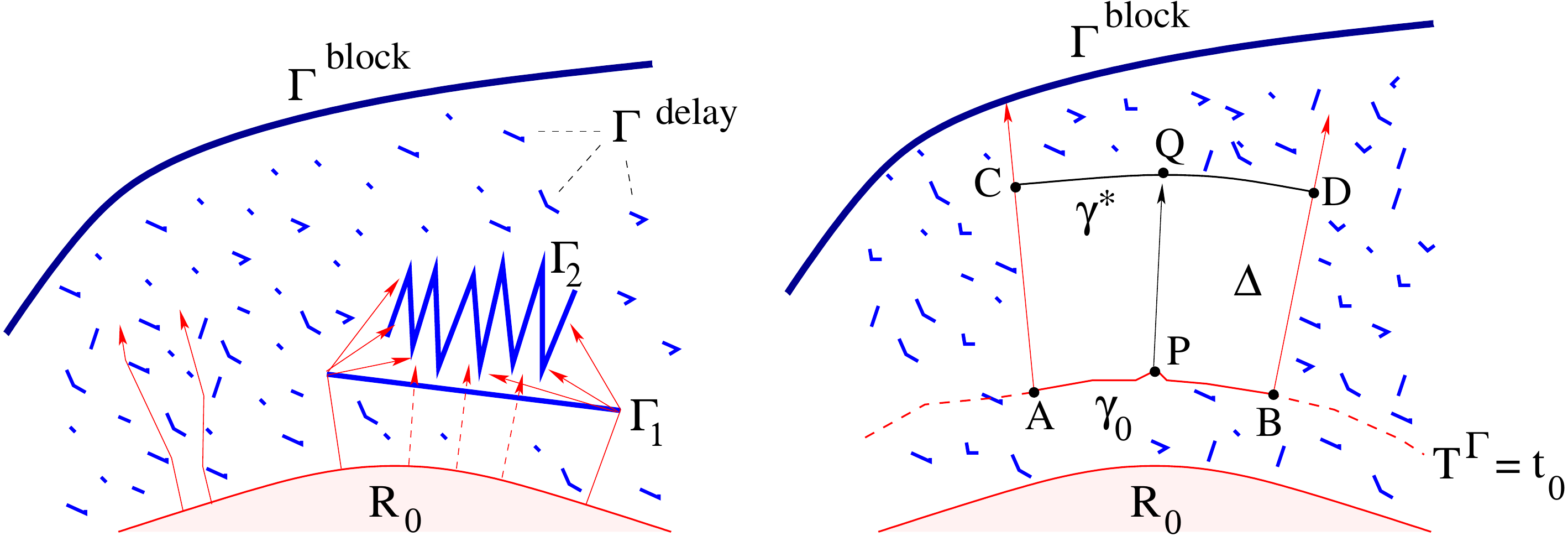}
   \caption{\small  Left: if the connected component $\Gamma_1$ of the barrier is removed, then there is not enough time to construct $\Gamma_2$, before the fire reaches it.  Right: by removing all barriers inside 
   a carefully chosen ``flow box"
   $\Delta$, the remaining portion $\Gamma^\diams= \Gamma\setminus\Delta$ still 
   form an admissible barrier, blocking  the fire within the same region as before
   and yielding a smaller total cost. }
   \label{f:fc181}
\end{figure}

This result is motivated by the following considerations. 
As shown in Fig.~\ref{f:fc181}, left,
the optimal barrier can be split as 
$$\Gamma~=~\Gamma^{\rm block}\cup \Gamma^{\rm delay}.$$
Here $\Gamma^{\rm block}= \partial \ov{R^\Gamma_\infty}$ is the portion  
which actually separates the burned region from the unburned one.
On the other hand, $\Gamma^{\rm delay}$ accounts for the walls whose only purpose is to delay the advancement of the fire front. 
Eventually, these walls are encircled by the fire on both sides.

We recall that the fire propagates with speed 1, while the barrier is constructed  
at speed $\sigma>1$.
Building a connected component $\Gamma_1$ of the barrier, with length $\ell_1$,
thus requires an amount of time $\tau_1 = \ell_1/\sigma$.
On the other hand, the fire needs up to time $\ell_1$ in order to 
completely surround $\Gamma_1$.    In some cases, it can 
thus be an advantage to 
construct some barriers with the sole purpose of slowing down the propagation of the fire.

At an intuitive level, however, 
building a barrier which contains a large number
of very small connected components should be ineffective, because the fire can
quickly get around each connected portion.   
To prove Theorem~\ref{t:1}, we need to show that a collection of walls which is 
dense on an open set cannot be optimal.  Indeed, some of these walls
should be removed, because   the time needed to build them is longer
 than the amount by which they delay the advancement of the fire front.
 
The heart of the matter is to understand which portions of the barrier can be removed.
As shown in Fig.~\ref{f:fc181}, left, the connected component 
$\Gamma_1$ delays the advancement of the fire front. 
If we remove $\Gamma_1$, then we do not have enough  time to construct $\Gamma_2$.
Hence, to achieve an admissible barrier $\Gamma'\subseteq\Gamma\setminus \Gamma_1$ satisfying (\ref{ad1}), 
we should  also remove the component $\Gamma_2$.   In turn, this may force us to 
remove further components $\Gamma_3, \Gamma_4,\ldots ~$
If at the end of this process we need to remove the outer component $\Gamma^{\rm block}$ as well, then the entire construction fails.

Toward a proof of Theorem~\ref{t:1}, 
we shall construct a ``flow box" $\Delta$, as shown in 
Fig.~\ref{f:fc181}, right.  Here the  lower boundary coincides with the location of the fire
front  at some time $t_0>0$. The two sides are straight lines, consisting
of optimal trajectories for the fire which do not intersect any of the barriers. The upper
boundary is a curve $\gamma^*$, consisting of points having a fixed distance $h>0$ from $\gamma_0$. These are the points that the fire would reach at 
time $t^*=t_0+h$, if no barriers were present.
A careful analysis will show that, by removing all the barriers contained inside $\Delta$, 
the remaining portion
$\Gamma^\diams~=~\Gamma\setminus\Delta$ is still admissible, and achieves a lower
total  cost (\ref{cost}).
\v
The remainder of the paper is organized as follows.  
Section~\ref{s:2} is concerned with 
the minimum time problem for the fire, in the presence of barriers.  
The first main result, Lemma~\ref{l:22}, considers a path $\xi$ that
crosses some of the connected components $\Gamma_i$ of the barrier.    By inserting additional loops, we prove the existence of a modified path
$\tilde\xi$ which starts and ends at almost the same points as $\xi$,
and does not touch the barrier.  Moreover, the difference between the lengths of the two paths is no greater than the total length $\sum_i m_1(\Gamma_i)$
of the components which were crossed.
The second main result of this section, Lemma~\ref{l:29}, shows that 
the set of times where the fire front touches a component $\Gamma_i$
is contained in an interval $[a_i, b_i]$ of length $b_i-a_i\le m_1(\Gamma_i)$.
Moreover, when no barrier is touched, the set reached by the fire expands with unit speed in all directions.   All these results are intuitively obvious when 
$\Gamma$ contains finitely many compact, connected components.
However, if $\Gamma$ is the union of countably many components,  possibly everywhere dense, a more careful proof is needed.

In Section~\ref{s:3} we prove some lemmas describing how the 
minimum time function $T^\Gamma$ in (\ref{TG}) 
changes when the barrier $\Gamma$ is perturbed.
This analysis is useful, because it allows us to approximate an arbitrary barrier 
with a polygonal one.

Section~\ref{s:4} continues the study of optimal trajectories for the fire, 
reaching points $x\in \R^2$ in minimum time without crossing 
the barrier $\Gamma$.  
The key result in this section (Lemma~\ref{l:5}) shows that, 
if the total length of all barriers is small, most
of these optimal trajectories for the fire contain long straight segments.
This fact can be rigorously stated in terms of an integral inequality.
The proof is first achieved in the case of polygonal barriers.   
The general case follows by 
an approximation argument.

Section~\ref{s:5} contains another key estimate. Roughly speaking,
Lemma~\ref{l:45} shows that, if a barrier $\Gamma$
is ``$\ve$-sparse", then the additional time needed by the fire to go around it
is bounded by $9\ve\, m_1(\Gamma)$.   We observe that the time needed
to construct this barrier is $\sigma^{-1} m_1(\Gamma)$, where $\sigma$ 
is the construction speed.    If $\ve>0$ is sufficiently small,
the time needed to construct this portion of barrier
is not compensated by its effectiveness in delaying the advance of the fire front.
One can thus conclude that  the barrier is not optimal.

The proof of Theorem~\ref{t:1} is then completed in  Section~\ref{s:6}.
It consists of two main steps.  First, we use
Lemma~\ref{l:5} to construct a ``flow box"
$\Delta$, as shown in Fig.~\ref{f:fc181}, whose sides are 
segments contained in optimal trajectories for the fire which do not cross the 
barrier $\Gamma$.   We then use Lemma~\ref{l:45} and show that, by removing all the portions of the barrier contained  inside $\Delta$, one obtains a new
admissible barrier $\Gamma^\diams= \Gamma\setminus \Delta$,  which yields
a smaller total cost.

\section{Optimal trajectories for the fire}
\label{s:2}
\setcounter{equation}{0}

\begin{figure}[htbp]
   \centering
 \includegraphics[width=0.24\textwidth]{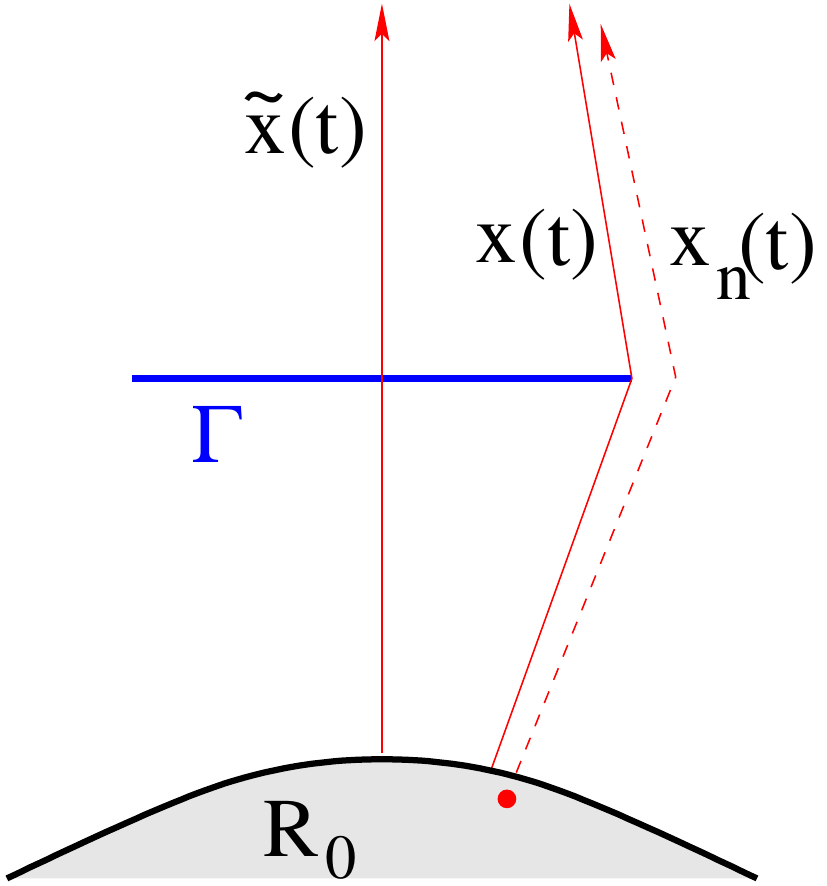}
   \caption{\small The trajectory $t\mapsto x(t)$ touches the barrier $\Gamma$, but does not cross it.
 Indeed, it can be obtained as a uniform limit of trajectories $x_n(\cdot)$ 
   that do not touch $\Gamma$. 
On the other hand, the trajectory $\tilde x(\cdot)$ is not admissible:
it goes right across the barrier. 
 }
   \label{f:fc149}
\end{figure}

In this section we focus on the optimization problem for the fire.  
Let $R_0\subset\R^2$ be a bounded, connected open set, and 
let $\Gamma=\cup_i\Gamma_i$ be a barrier, consisting of countably many
compact, rectifiable, connected components, with finite total length. We seek trajectories that, starting from 
the closure $\ov {R_0}$, reach points $x\in \R^2$ in minimum time, without
crossing $\Gamma$.   To achieve the 
existence
of these optimal trajectories, referring to Fig.~\ref{f:fc149} we introduce 
\begin{definition}\label{d:adm}
A trajectory for the fire $t\mapsto x(t)$, $t\in [0, T]$, is {\bf admissible}
if there exists a sequence of 1-Lipschitz trajectories $t\mapsto x_n(t)$ such that 
\bel{aco1}
x_n(0)\in R_0\,,\qquad x_n(t)\notin \Gamma\qquad\forall t\in [0,T],
\eeq
and moreover $x_n(t)\to x(t)$ uniformly on $[0,T]$, as $n\to\infty$.

We say that a trajectory $t\mapsto x(t)$ 
{\bf does not touch the barrier} $\Gamma$ if $x(t)\notin \Gamma$ 
for all times $t\geq 0$.
If $x(\cdot)$ is the uniform limit of trajectories $x_n(\cdot)$ that do not touch
$\Gamma$, we say that $x(\cdot)$ {\bf does not cross the barrier} $\Gamma$.
\end{definition}

Given a point $\bar x\in \R^2$, we seek an admissible trajectory $t\mapsto x(t)$ which starts from a point in the closure $\ov{R_0}$ and
reaches  $\bar x$ in minimum time without crossing $\Gamma$. 
If $\bar x$ can be reached in 
finite time, the existence of such an optimal trajectory is 
straightforward. 
Indeed, define
$$\bega{rl}T_{\rm inf}&\doteq~\ds\lim_{\ve\to 0} \Big[\hbox{infimum time needed to reach a point in the ball $B(\bar x, \ve)$,}\\[3mm]
&\qquad\qquad\qquad 
 \hbox{ starting from $R_0$ and without touching the barrier $\Gamma$}\Big].\enda$$
Let $x_n:[0, T_n]\mapsto \R^2$ be a minimizing sequence of 1-Lipschitz trajectories, 
satisfying (\ref{aco1}) together with 
$$x_n(T_n)~\to ~\bar x,\qquad\qquad T_n~\to~T_{\rm inf}\qquad\qquad 
\hbox{as}\quad n\to\infty.$$
By taking a subsequence we can assume the uniform convergence 
$x_n\to x$, for some limit function $x:[0, T_{\rm inf}]\mapsto \R^2$. 
According to Definition~\ref{d:adm}, this limit trajectory is admissible.  Hence it provides an optimal solution.
\v
Given a trajectory $\xi:[0,\tau]\mapsto\R^2$ that crosses part of the barrier,
the next lemma provides the key tool for constructing trajectories that
``loop around" each connected component, and reach almost the same 
endpoint without touching $\Gamma$. 

\begin{lemma}
\label{l:22} 
Consider a barrier $\Gamma= \cup_{i\geq 1}\Gamma_i$, written as the union of 
its connected components.  Assume that $\R^2\setminus \Gamma$ is connected.
Let $\xi:[0,\tau]\mapsto\R^2$ be a Lipschitz path, parameterized by arc length, such that
\bel{xgi}
\xi(t)\notin \Gamma_i\qquad\forall t\in [0,\tau],~~i\leq\nu.\eeq
Then, for any $\epsilon>0$, there exists a path $\Tilde\xi:[0,\Tilde\tau]\mapsto\R^2$,  also  parameterized by arc length, such that 
\bel{tx1}\bigl|\Tilde\xi(0)-\xi(0)\bigr|~\leq~\epsilon,\qquad \bigl|\Tilde\xi(\Tilde\tau)-\xi(\tau)\bigr|~\leq~\epsilon,\eeq
\bel{tx2}\Tilde\xi(t)\notin\Gamma\qquad\forall t\in [0,\tilde \tau],\eeq
and with length
\bel{tx4}\Tilde\tau~\leq~\tau+ \sum_{i>\nu} m_1(\Gamma_i)\,.\eeq
\end{lemma}
\v

\begin{figure}[htbp]
   \centering
 \includegraphics[width=0.6\textwidth]{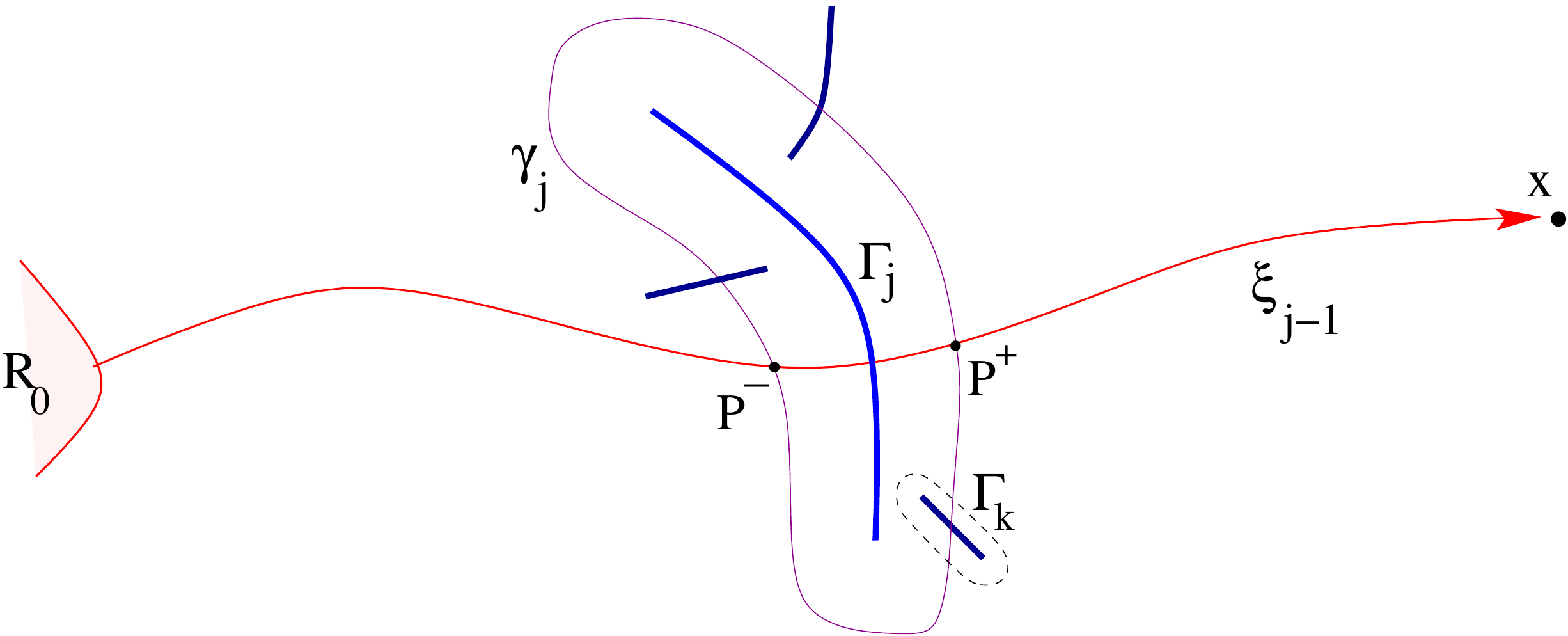}
   \caption{\small The construction used in the proof of Lemma~\ref{l:27}.
   If the trajectory $\xi_{j-1}(\cdot)$ crosses the set $\Gamma_j$,
   we construct a detour $\gamma_j$ of radius $r_j$ around $\Gamma_j$.
   At a subsequent step, we may be forced to construct a second detour
to avoid hitting the component $\Gamma_k$.   Hence the new path may get closer to $\Gamma_j$. In the inductive construction, it is essential 
to show that all paths keep a uniformly positive distance from $\Gamma_j$. }\label{f:fc186}
\end{figure}

{\bf Proof.}  Let $\epsilon>0$ be given.
The new path $\Tilde \xi$ will be obtained as limit
of a sequence of paths $\xi_j:[0,\tau_j]\mapsto\R^2$, $j\geq 0$, by an inductive procedure.
Each inductive step will also determine two auxiliary constants $r_j, \delta_j>0$.
\v
{\bf 1.} 
The induction starts by setting 
$\tau_0=\tau$, and defining $\xi_0(t)=\xi(t)$ for all $t\ge 0$. Moreover, we choose
$\delta_0>0$ so that
\bel{del0}
\delta_0\,<\, {\epsilon\over 4}\,,\qquad \delta_0~<~d(\xi(t),\,\Gamma_i)\quad \forall
t\in [0,\tau], ~~i=1,\ldots,\nu.\eeq
\v
For every $j\geq 1$, the constants $r_j, \delta_j>0$ and the path $\xi_j:[0,\tau_j]\mapsto\R^2$ will satisfy
the following properties.
\begi
\item[(i)] For every $i\leq j$ and $t\in [0,\tau_j]$ one has
\bel{dg1} d\bigl(\xi_j(t), \Gamma_i\bigr)\,\geq\, (2-2^{i-j})\delta_i\,.\eeq
\item[(ii)] For $j\leq\nu$ we simply take $\xi_j=\xi_0$.  
For $j>\nu$, the length of the path $\xi_j$ satisfies
\bel{dg2}
\tau_j~<~\tau_{j-1}+(1+\epsilon) m_1(\Gamma_j).\eeq
\item[(iii)] The endpoints satisfy
\bel{dg3} |\xi_j(0)-\xi(0)|\,\leq \,(1-2^{-j})\epsilon,\qquad\quad  |\xi_j(\tau_j)-\xi(\tau)|\,\leq \,(1-2^{-j})\epsilon.\eeq
\item[(iv)] The constant $r_j$ is chosen so that
 $r_j<\delta_{j-1}/4$.  Moreover, every component $\Gamma_k$
which intersects $B(\Gamma_j, 2r_j)$  has length $m_1(\Gamma_k)<\delta_{j-1}/4$.
\item[(v)] The constant
$\delta_j\in \,]0,r_j/4]$ is chosen so that, for every $k\geq 1$ such that  
\bel{dg5}
\Gamma_k\cap \Big( \ov B(\Gamma_j, 2r_j)\setminus B(\Gamma_j, r_j/2)
\Big)\,\not=\,\emptyset,\eeq 
one has
\bel{dg4} 4\delta_j~\leq~
d(\Gamma_k,\Gamma_j)~\doteq~
\min\,\Big\{|x-y|\,;~x\in \Gamma_j, ~y\in \Gamma_k
\Big\} .\eeq
\endi

Note that, even if we choose $\xi_j=\xi_0$ for $j=1,\ldots,\nu$, it is not possible
to start the induction procedure at $j=\nu$. Indeed, the initial steps must be
performed in order to define suitable constants $r_j, \delta_j$, $j=1,\ldots,\nu$.
\v
{\bf 2.} Assuming that the induction has been completed up to step $j-1$, we 
describe how to accomplish step $j$.

Consider the   path $\xi_{j-1}:[0,\tau_{j-1}]$, and  the connected component
$\Gamma_j$.   For a given radius $r>0$, define
$$\gamma_j~\doteq~[\hbox{boundary of the unbounded connected component
of $\R^2\setminus B(\Gamma_j, r)$}].$$
This is a simple closed curve, that winds around $\Gamma_j$, and has length $\approx 2 m_1(\Gamma_j)$.
We choose $r=r_j>0$  small enough, so that the following holds:
\bel{me1}m_1(\gamma_j)~<~(2+\epsilon)m_1(\Gamma_j),\eeq
\bel{me0}r_j~<~\ds {\delta_{j-1}\over 2}\,,\qquad\qquad r_j~<~{1\over 4}
\,\min_{1\leq i<j} d(\Gamma_i, \Gamma_j),\eeq
and moreover 
\begi
\item[{\bf (P$_j$)}] {\it 
Every connected component $\Gamma_k$, $k\not= j$, that intersects 
$\ov B(\Gamma_j, 2r_j)$ has length 
$<\delta_{j-1}/4$.}

\item[{\bf (P$_j'$)}] {\it Every disc of radius $\delta_{j-1}$ intersects the 
unbounded connected component of $\R^2\setminus B(\Gamma_j, r_j)$.}
\endi

Note that all the above can certainly be achieved, because there are 
only finitely many components
$\Gamma_k$ 
whose length is $> \delta_{j-1}/4$.   Choosing $r_j>0$ small enough, $\gamma_j$
will not intersect any of them.  Moreover, the  inequality (\ref{me1}) follows by
well known results in geometric measure theory \cite{AFP, Federer}.
Indeed, since $\Gamma_j$ is rectifiable, the neighborhoods of radius $r$ around
$\Gamma_j$ satisfy
$$\lim_{r\to 0}~m_2\left({B(\Gamma_j,r)\over 2 r}\right)~=~m_1(\Gamma_j).$$
Hence the co-area formula yields
\bel{ca}\liminf_{r\to 0} ~m_1\Big(\partial B(\Gamma_j,r)\Big)~\leq~2m_1(\Gamma).
\eeq
We then choose a constant $\delta_j$ according to (v) above.

Finally, the new path $\xi_j:[0, \tau_j]\mapsto \R^2$ is defined as follows.
If $\xi_{j-1}(t)\notin \ov B(\Gamma_j, r_j)$ for all $t\in [0, \tau_{j-1}]$, 
then there is no need to modify the previous path, and we can simply set
$$\tau_j\,=\,\tau_{j-1}\,,\qquad\qquad \xi_j(t)=\xi_{j-1}(t).$$

\begin{figure}[htbp]
   \centering
 \includegraphics[width=1.0\textwidth]{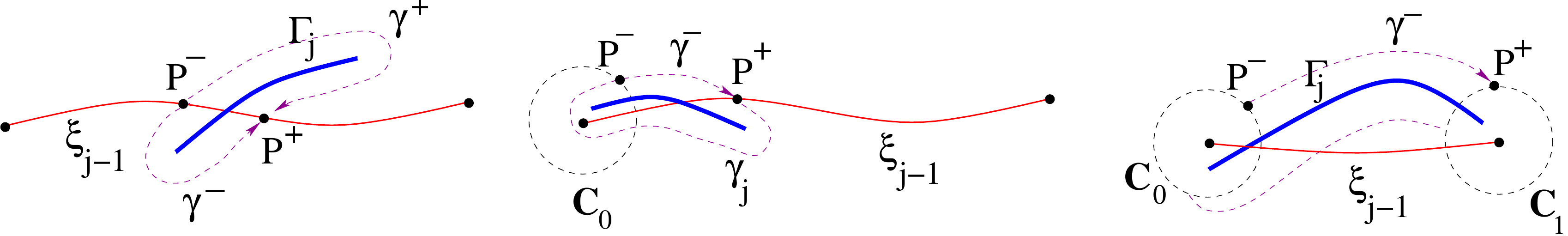}
   \caption{\small If the path $\xi_{j-1}$ intersects the component $\Gamma_j$ of the barrier, a detour must be constructed.    The figures on the left, center, and right illustrate Cases 1, 2, and 4, respectively.}\label{f:fc187}
\end{figure}

Otherwise, we add a detour so that the new path will remain bounded away from the 
component $\Gamma_j$ of the barrier. For this purpose,
define the times
$$\bega{l}
t^-~\doteq~\inf\,\bigl\{t\in [0,\tau_{j-1}]\,;~\xi_{j-1}(t)\in \ov B(\Gamma_j,r_j)
\bigr\},\\[3mm] t^+~\doteq~\sup\,\bigl\{t\in [0,\tau_{j-1}]\,;~\xi_{j-1}(t)\in \ov B(\Gamma_j,r_j)
\bigr\},\enda$$
and the points
$$P^-\,=\, \xi_{j-1}(t^-), \qquad P^+\,=\, \xi_i(t^+).$$
Various cases need to be considered (see Fig.~\ref{f:fc187}).

\v
CASE 1: $0<t^-\leq t^+<\tau_{j-1}$, shown in Fig.~\ref{f:fc187}, left.

In this basic case 
we observe that the two points $P^-, P^+$ 
divide the simple closed curve $\gamma_j$ into two parts, say $\gamma_j^-$, $\gamma_j^+$.  To fix the ideas, assume 
\bel{me2}s_j~\doteq~m_1(\gamma_j^-)~\leq ~m_1(\gamma_j^+).\eeq    
Let $s\mapsto \gamma_j^-(s)$ be an arc-length parameterization of 
$\gamma_j^-$, with 
$$\gamma_j^-(0)~=~ P^-,\qquad \qquad \gamma_j^-(s_j)~=~P^+.$$
We then define the new path $\xi_j$ by adding a detour around $\Gamma_j$ as follows:
\bel{gamj}
\xi_j(t)~=~\left\{
\bega{rl} \xi_{j-1}(t)\quad &\hbox{if}\quad t\in [0, t^-],\\[3mm]
\gamma_j^-(t-t^-)\quad &\hbox{if} \quad t\in [t^-,~t^-+ s_j],\\[3mm]
\xi_{j-1}(t - s_j+t^+-t^-)\quad &\hbox{if} \quad t\in [t^-+s_j, ~\tau_j].
\enda\right.
\eeq
Here 
$$\tau_j~=~\tau_{j-1} + s_j -(t^+-t^-) .$$
\v

CASE 2: $0=t^-\leq t^+<\tau_{j-1}$, shown in Fig.~\ref{f:fc187}, center.

Call $P^+=\xi_{j-1}(t^+)$.    Observe that, by (\ref{me0}) and 
{\bf (P$_j'$)}, the curve 
$\gamma_j$ has non-empty intersection with the circumference 
\bel{circ0}{\bf C}_0\doteq
\{ y\in\R^2\,;~|y-\xi_{j-1}(0)|=
\delta_{j-1}\}.\eeq    Therefore, starting at $P^+$ and moving along the simple closed curve $\gamma_j$, we can reach some point $P^-$ on ${\bf C_0}$ in two different ways: clockwise and counterclockwise.   We choose $\gamma^-\subset\gamma_j$
to be the shortest among these two paths.

As shown in Fig.~\ref{f:fc187}, center, we parameterize $\gamma^-$ by arc length, so that
$$\gamma^-(0)= P^-,\qquad \gamma^-(s_j)=P^+$$
for some $s_j>0$.
Choosing $t^-$ so that  $P^-=\xi_{j-1}(t^-)$,
we define the new path $\xi_j(\cdot)$ by setting $\tau_j= \tau_{j-1}- t^-+s_j$ and
\bel{xi4}
\xi_j(t)~=~\left\{\bega{cl} \gamma^-(t)\quad &\hbox{if}~~t\in [0, s_j],\\[3mm]
\xi_{j-1}(t+t^- - s_j)\quad&\hbox{if}~~t\in [s_j, ~\tau_{j-1}- t^-+s_j].\enda\right.
\eeq
\v
CASE 3: $0<t^-\leq t^+=\tau_{j-1}$.

 In this case, the new path $\xi_j(\cdot)$
will connect $\xi_{j-1}(0)$ with a point on the circumference
\bel{circ1} {\bf C_1}\doteq
\bigl\{ y\in\R^2\,;~|y-\xi_{j-1}(\tau_{j-1})|=
\delta_{j-1}\bigr\}.\eeq  Since this is entirely similar to Case 2, we omit the details.
\v
CASE 4: $0=t^-<t^+= \tau_{j-1}$, shown in Fig.~\ref{f:fc187}, right.

In this case, the simple closed curve $\gamma_j$ intersects both 
circumferences ${\bf C}_0$ and ${\bf C}_1$ in (\ref{circ0}), (\ref{circ1}).
We then choose a point $P^-\in {\bf C}_0$ and a point $P^+\in {\bf C}_1$
so that the portion $\gamma^-\subset \gamma_j$ connecting $P^-$ with
$P^+$ is as short as possible.   

We now parameterize $\gamma^-$ by arc length, so that
$$\gamma^-(0)= P^-,\qquad \gamma^-(s_j)=P^+$$
for some $s_j>0$.
The new path $\xi_j(\cdot)$ is defined simply by setting $\tau_j=s_j$ and
\bel{xi5}
\xi_j(t)~=~\gamma^-(t)\qquad \forall ~ t\in [0, s_j].
\eeq

\v
{\bf 3.} Having constructed a sequence of paths $\xi_j:[0,\tau_j]\mapsto\R^2$,
$j\geq 0$, by taking the limit as $j\to \infty$ we will obtain
a path $\Tilde \xi(\cdot)$ which satisfies
the properties (\ref{tx1})-(\ref{tx2}), together with 
\bel{tx3}\Tilde\tau~\leq~\tau+ (1+\epsilon)\sum_{i>\nu} m_1(\Gamma_i)\,.\eeq
  For convenience, we extend the definition
of each path $\xi_j$ to all of $\R_+$  by setting 
$$\xi_j(t)\,=\, \xi_j(\tau_j)\qquad\quad\forall t\geq \tau_j\,.$$
Toward a proof of  (\ref{tx1}), we observe that our construction implies
$$|\xi_j(0)-\xi_{j-1}(0)\bigr|~\leq~\delta_{j-1}~<~2^{-j}\epsilon\,.$$
$$|\xi_j(\tau_j)-\xi_{j-1}(\tau_{j-1})\bigr|~\leq~\delta_{j-1}~<~2^{-j}\epsilon\,.$$
Recalling that $\xi_0=\xi$ and  summing these inequalities from $1$ to $j$, we obtain (\ref{dg3}).
\v
{\bf 4.}
Our construction guarantees that, for  every $j\geq 1$, the length $\tau_j$ of the new curve $\xi_j(\cdot)$ satisfies (\ref{dg2}).  Notice that, if Case 4 occurs,
we have the even sharper bound 
$$\tau_j~\leq~(1+\epsilon) m_1(\Gamma_j).$$
In addition, since we are assuming that the initial path $\xi(\cdot)$
does not intersect any of the components $\Gamma_1,\ldots, \Gamma_\nu$,
our choice of  $\delta_0>0$ in  (\ref{del0}) guarantees that 
no modification need to be 
done in the first $\nu$ steps of the algorithm.
Hence 
$\xi_j(\cdot)=\xi_{j-1}(\cdot)$ for all $j=1,\ldots,\nu$.   In particular, this implies
\bel{tnu}\tau_j~=~\tau\qquad\forall j=1,\ldots,\nu.\eeq

If one of the Cases 1-2-3 occurs, then 
our construction yields
\bel{xijj}
\bigl| \xi_j(t) - \xi_{j-1}(t)|~\leq~(1+\epsilon)\,m_1(\Gamma_j)\qquad\forall t\geq 0.\eeq
In Case 4 the above estimate can fail.
However,  by (\ref{dg3}), Case 4 in the above construction can occur
only finitely many times.   Indeed, there are at most finitely many
connected components $\Gamma_j$ of length 
$$m_1(\Gamma_j)~\geq~|\xi(\tau)-\xi(0)| - 4\epsilon.$$
We thus conclude that the sequence $\xi_j(\cdot)$ is Cauchy. 
As $j\to\infty$, we have the  convergence
$\xi_j(t)\to  \xi_\infty(t)$, uniformly for $t\geq 0$.
\v
{\bf 5.} Since all paths $\xi_j$ are 
Lipschitz continuous with constant 1, the limit  path $\xi_\infty$ is  1-Lipschitz
as well.    We can now parameterize $ \xi_\infty$
 by arc-length, and obtain a path $\Tilde\xi:[0,\Tilde\tau]\mapsto\R^2$,
with 
$$\Tilde\tau~\leq ~\liminf_{j\to\infty} \,\tau_j~\leq~\tau +(1+\epsilon) \sum_{j>\nu} m_1(\Gamma_j).$$
Notice that the above estimates follow from (\ref{tnu}) and (\ref{dg2}).
This proves (\ref{tx3}).   

The bounds (\ref{tx1}) are an immediate consequence of (\ref{dg3}).
\v
{\bf 6.} In this step we show that (\ref{tx2}) holds.   Namely,  
the path $\Tilde\xi$ does not touch any of the components $\Gamma_i$, 
$i\geq 1$ of the barrier.

This claim will be proved by showing that, 
for a fixed $j\geq 1$ and every $k\geq j$
one has
\bel{dkj}
d(\xi_k(t), \Gamma_j)~\doteq~\min\bigl\{ |\xi_k(t)-y|\,;~~y\in \Gamma_j\bigr\}~\geq ~\delta_j\qquad\quad\forall t\in [0, \tau_k]\,.\eeq

By construction, it immediately follows that 
$$d(\xi_j(t), \Gamma_j)~\geq~r_j~\geq~4\delta_j\qquad\forall t\in [0,\tau_j].$$

We now observe that, for
 any $k>j$, the path $\xi_k(\cdot)$ is obtained from $\xi_j(\cdot)$
by replacing some of its sections by 
detours  $\gamma_i^-(\cdot)$,~ $i=j+1,\ldots, k$, ~where
$$d(\gamma_i(s), \Gamma_i)~\doteq~\min\,\Big\{
\bigl|\gamma(s)-y\bigr|\,;~~y\in\Gamma_i\Big\}~=~r_i\qquad\forall s.$$
Three cases must be considered.

CASE 1: (\ref{dg5}) holds, and hence by construction (\ref{dg4}) holds as well.
In this case we have
$$d(\gamma_k(s), \Gamma_j)~\geq~ d(\Gamma_k,\Gamma_j) - d(\gamma_k(s),
\Gamma_k)~\geq~4\delta_j - r_k~\geq~3\delta_j\,.$$
\v
CASE 2: $d(\Gamma_k,\Gamma_i)>2r_j$.   In this case, for every $s$  
we trivially have
$$d(\gamma_k(s), \Gamma_j)~\geq~
d(\Gamma_k,\Gamma_j) - d(\gamma_k(s),
\Gamma_k)~\geq~2r_j-r_k~>~r_j~\geq ~4\delta_j\,.$$

CASE 3: $\Gamma_k\subseteq \ov B(\Gamma_j, r_j/2)$.
If $\xi_k(\cdot)=\xi_{k-1}(\cdot)$, the conclusion (\ref{dkj}) follows
by induction on $k$.    It thus suffices to consider the case where
$\xi_k(\cdot)$ is obtained from $\xi_{k-1}(\cdot)$ by inserting some 
nontrivial portion
of a curve $\gamma_k\subseteq\{x\,;~d(x,\Gamma_k)=r_k\}$.

In this case, our algorithm implies that there exists a finite
sequence
$$j = i(0)<i(1)<\cdots <i(N)= k$$
such that every 
curve $\gamma_{i(\ell)}~\doteq~\{x\,;~d(x,\Gamma_k)=r_k\}$
intersects the previous one:
$$\gamma_{i(\ell)}\cap \gamma_{i(\ell-1)}~\not= ~\emptyset\qquad\forall 
\ell=1,\ldots,N.$$
Considering the diameters of the sets $\gamma_{i(\ell)}$, we thus have the bound
\bel{gij}\bega{rl}\ds
\min_s d\bigl(\xi_k(s), \Gamma_j\bigr)&\ds\geq~
\min_s d\bigl(\gamma_{i(1)}(s), \Gamma_j\bigr) - \sum_{\ell=2}^N
\hbox{diam}(\gamma_{i(\ell)})\\[4mm]
&\ds\geq~\bigl(d(\Gamma_{i(1)},\Gamma_j) - r_{i(1)}\bigr) - \sum_{\ell=2}^N
\bigl(2 r_{i(\ell)}+ m_1(\Gamma_{i(\ell)})\bigr).\\[4mm]
&\ds\geq~\left(4\delta_j - {\delta_j\over 2}\right) - \sum_{\ell=2}^N\left(
2\cdot  2^{j-i(\ell)} {\delta_j\over 4} + 2^{j-i(\ell)} {\delta_j\over 4}\right)~>~\delta_j\,.
\enda\eeq

Indeed, the condition (\ref{dg5}) applies to $\Gamma_{i(1)}$, hence 
(\ref{dg4}) holds.  Moreover, using the property
{\bf (P$_j$)}  with $j$ replaced by $i(1), \ldots, i(N)$, we obtain
$$m_1(\Gamma_{i(\ell)})~\leq~{\delta_{i(\ell)-1}\over 4}~\leq~2^{j-i(\ell)} \,{\delta_j\over 4}\,.$$

Combining the above three cases, we conclude that (\ref{dkj}) holds. 
Taking the limit as $k\to\infty$, we conclude that $d(\Tilde\xi(t), \Gamma_j)\geq \delta_j$
for all $t\geq 0$ and $j\geq 1$. This establishes  (\ref{tx2}), 
\v
{\bf 7.}
To obtain the bound (\ref{tx4}) on the length of the new path,
define
$$\Hat \tau~\doteq~\min\left\{ \Tilde\tau,~\tau +\sum_{i>\nu} m_1(\Gamma_i)\right\}.
$$
By (\ref{tx3}), we trivially have
$$|\Tilde \tau -\Hat\tau|~\le~\epsilon \sum_{i>\nu} m_1(\Gamma_i).$$
Therefore, if we replace the path $\tilde \xi$ by its restriction to the subinterval 
$[0, \Hat\tau]$, the conditions (\ref{tx2})-(\ref{tx4}) are  satisfied, while
the second inequality in (\ref{tx1}) will be replaced by 
$$\bigl|\Tilde\xi(\Hat \tau)- \xi(\tau)\bigr|~\leq~\epsilon + \epsilon\sum_{i>\nu} m_1(\Gamma_i) ~\leq~\epsilon (1+m_1(\Gamma)).$$
Since $\epsilon>0$ can be chosen arbitrarily small, this completes the proof.
\endproof
\v

The next result is concerned with  the length of the portion of the barrier
which is touched by the fire at a given time $t$.   We recall that, if $\Gamma$
is admissible, the linear bound  (\ref{ad1}) must hold.
\begin{lemma}
\label{l:21} Given  an admissible barrier $\Gamma$,
consider the function
\bel{vpdef}\vp(t)~=~m_1\bigl(\ov {R^\Gamma(t)}\cap\Gamma\bigr).\eeq
Then
\begi
\item[(i)] $\vp$ is nondecreasing and right continuous.
\item[(ii)] The set of times where the constraint is not saturated
\bel{nsat}
\U~\doteq~\Big\{t>0\,;~~m_1\bigl(\ov {R^\Gamma(t)}\cap\Gamma\bigr)\,<\,\sigma t\Big\}\eeq
is open.
\item[(iii)] The set of times where the constraint is saturated
\bel{sat}
\S~\doteq~\Big\{t\geq 
0\,;~~m_1\bigl(\ov {R^\Gamma(t)}\cap\Gamma\bigr)\,=\,\sigma t\Big\}\eeq
is closed.
\endi
\end{lemma}
\v
{\bf Proof.}
{\bf 1.} 
For any $0< t_1<_2$ we  have $R^\Gamma(t_1)\subseteq R^\Gamma(t_2)$.
Hence $\vp$ is nondecreasing.  
\v
{\bf 2.} 
Next, we claim that $\vp$ it is right continuous.  
Indeed, consider a decreasing sequence of times $t_n\downarrow t_0$.
Since the fire propagates with unit speed, we have
$$\ov {R^\Gamma(t_n)}~\subseteq~ \ov{B\Big(R^\Gamma(t_0),~t_n-t_0\Big)}\,$$
hence
$$\ov {R^\Gamma(t_0)}\cap\Gamma~=~\bigcap_{n\geq 1} \bigl(
\ov {R^\Gamma(t_n)}\cap\Gamma\bigr).
$$
The right continuity of $\vp$ now follows from the dominated convergence theorem.
\v
{\bf 3.} By the previous two steps it follows that $\vp$ is upper semicontinuous.
Hence the function $t\mapsto \vp(t)-\sigma t$ is upper semicontinuous as well.
We thus conclude that the set $\U$ where $\vp(t)-\sigma t <0$ is open.
The closure of $\S= \R_+\setminus\U$ follows immediately.
\endproof
\v

The next lemma will play a key role in the sequel. 
The intuitive idea is simple: let $t=a_i$ be the first time when the fire front touches
the connected component $\Gamma_i$.  Immediately afterwards, the fire 
starts going around
$\Gamma_i$, clockwise as well as counterclockwise, until this connected component is
completely surrounded.  This will happen at some time $b_i$ with $b_i-a_i\leq m_1(\Gamma_i)$. On the other hand, 
when the fire front does not touch any of the 
barriers $\Gamma_j$, it expands freely 
with unit speed in all directions. Therefore, the distance between level sets
of the time function $T^\Gamma$ increases at unit rate.  

\begin{lemma}\label{l:29} Consider a barrier $\Gamma= \cup_{i\geq 1}\Gamma_i$, 
written as the union of 
its compact connected components.  
Assume that $\R^2\setminus \Gamma$ is connected.
\begi
\item[(i)] For each $i\geq 1$, the set of times when the fire front touches $\Gamma_i$
\bel{around}
J_i~\doteq~\Big\{t\geq 0\,;~\partial\,\ov{R^\Gamma(t)}\cap \Gamma_i~\not= ~\emptyset
\Big\}\eeq
is contained within an  interval $[a_i, b_i]$  of length
$b_i-a_i\leq m_1(\Gamma_i)$.
\item[(ii)] For any $0\leq \tau<\tau'$, one has
\bel{lsb}B\bigl( R^\Gamma(\tau), r\bigr)~\subseteq~\ov{R^\Gamma(\tau')},
\qquad \hbox{with} \qquad r~=~
m_1\left( [\tau, \tau']\setminus \bigcup_{i\geq 1}[a_i, b_i]\right).\eeq
\endi
\end{lemma}
\v
{\bf Proof.} {\bf 1.} By the assumptions, each $\Gamma_i$ is simply connected.
Let
\bel{aidef} a_i~\doteq~\inf\,\Big\{ t\geq 0\,;~\ov{R^\Gamma(t)}\cap \Gamma_i \not= \emptyset\
\Big\}~=~\min_{x\in \Gamma_i} ~T^\Gamma(x)\eeq
be the first time when the fire touches $\Gamma_i$. By the lower semicontinuity of $T^\Gamma$  and the compactness of $\Gamma_i$, it is clear that $a_i$ is actually a minimum.   
We will prove part (i) of the lemma  by showing that,
at any time $\tau> a_i+m_1(\Gamma_i)$, the component $\Gamma_i$
is entirely contained in the interior of the set $\ov{R^\Gamma(\tau)}$.
\v
{\bf 2.} Toward our goal, we first choose  $\ve>0$  such that
\bel{epc}
4\ve~<~\tau - a_i -m_1(\Gamma_i),\eeq
then we choose an integer $\nu>i$ so large that 
\bel{lsm}
\sum_{k>\nu}m_1(\Gamma_k)~<~\ve \,.\eeq
Finally, we choose a radius $0<\rho<\ve$ small enough so that
\bel{nint}
B(\Gamma_i,\rho)\,\cap\,\Gamma_j~=~\emptyset\qquad 
\forall~ j=1,\ldots,\nu, ~~j\not= i.\eeq
With the above choices, we will show that
\bel{bgr}B(\Gamma_i, \rho)~\subset~\ov{R^\Gamma(\tau)}.\eeq
\v
{\bf 3.}   To prove (\ref{bgr}), fix any point $x\in B(\Gamma_i, \rho)\setminus \Gamma_i$.
For $0<r<\rho$,
consider the open neighborhood 
$B(\Gamma_i, r)$ of radius $r$ around $\Gamma_i$.
By a suitable choice of $r>0$, we claim that the following properties can be achieved.
\begi
\item[(i)] Calling $\gamma$ the boundary of the unbounded connected component of 
$\R^2\setminus B(\Gamma_i, r)$, we have
\bel{mg1}
m_1(\gamma)~<~2m_1(\Gamma_i) + \ve\, .\eeq
 \item[(ii)] The point $x$ lies in the unbounded connected component of 
$\R^2\setminus B(\Gamma_i, r)$.
\endi
Indeed, the property (i) follows by the same argument used in (\ref{ca}).
The property (ii) follows from the fact that $\Gamma_i$ is compact and simply connected, while 
$x\notin \Gamma_i$.

\begin{figure}[htbp]
   \centering
 \includegraphics[width=0.45\textwidth]{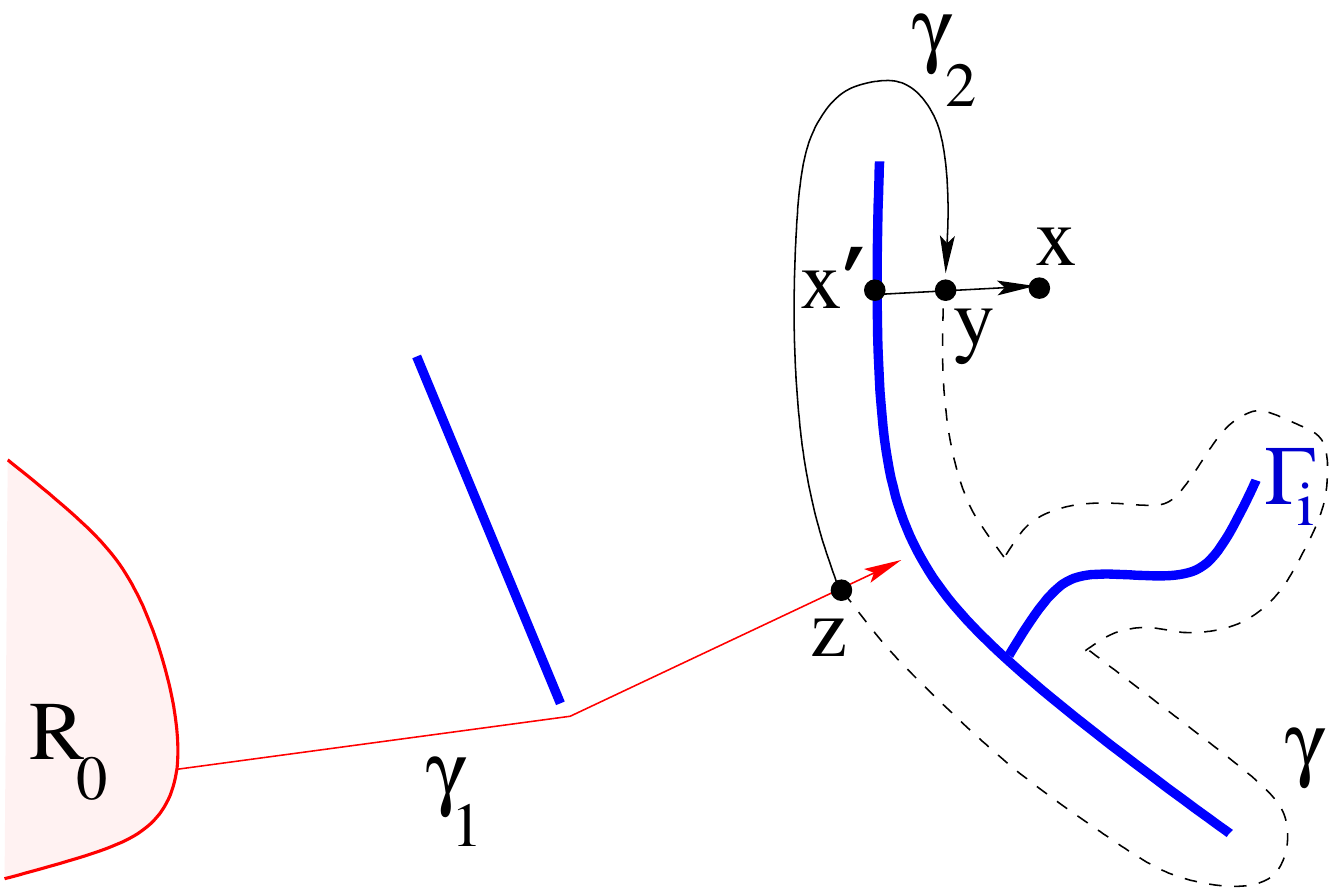}
   \caption{\small The construction used in the proof of part (i) of Lemma~\ref{l:29}.}
   \label{f:fc196}
   \end{figure}

\v
{\bf 4.} As shown in Fig.~\ref{f:fc196}, let $x'\in \Gamma_i$ be one of the points closest to $x$,
so that $|x-x'|<\rho$.   By construction, the segment with endpoints $x',x$ intersects the 
simple closed curve $\gamma$ at least at one point, say $y\in \gamma$.
This implies
\bel{xyd}
|x-y|~<~|x-x'|~<~\rho~\leq ~\ve.\eeq
Next, since $\ov{R^\Gamma(a_i)}\cap\Gamma\not=\emptyset$, 
there exists a trajectory for the fire that starts inside $R_0$ and crosses the curve $\gamma$
at some point $z$ before time $t=a_i$.

We now consider the path $\xi:[0,\ell]\mapsto\R^2$ obtained by concatenating the following three paths:
\begi
\item A path $\gamma_1$, starting inside $R_0$ and reaching $z\in \gamma$ without crossing
the barrier $\Gamma$.   This path has length $\ell_1< a_i$.
\item A path $\gamma_2$ contained within $\gamma$, starting at $z$ and ending at $y$.
Since we can move along $\gamma$ both clockwise or counterclockwise, by choosing the
shorter path we can assume that $\gamma_2$ has length 
$$\ell_2~\leq~ {1\over 2} m_1(\gamma)
~<~ m_1(\Gamma_i) + \ve.$$
\item A path $\gamma_3$ consisting of the segment with endpoints $y,x$.
By (\ref{xyd}), its length is $\ell_3 < \ve$.
\endi
The total length of this path $\xi(\cdot)$ is thus
$$\ell~ =~ \ell_1+\ell_2+\ell_3~<~a_i + \bigl(m_1(\Gamma_i) +\ve\bigr) +\ve~<~\tau -2\ve.$$
Notice that, by construction, the path $\xi$ does not cross any of the components $\Gamma_1,\ldots,\Gamma_\nu$.    Applying Lemma~\ref{l:22}, for any $\ve'>0$ we can find a new path
$\Tilde\xi:[0,\Tilde \ell]\mapsto\R^2\setminus\Gamma$ such that
$$\Tilde\xi(0)\in R_0\,,\qquad\qquad \bigl|\Tilde\xi(\Tilde \ell)-x\bigr|~<~\ve',$$
and moreover
$$\Tilde \ell~\leq~a_i + m_1(\Gamma_i) +2\ve + \sum_{k>\nu} m_1(\Gamma_k) ~<~\tau.$$
This implies $x\in \ov{R^\Gamma(\tau)}$, as claimed.   Hence part (i) of the lemma is proved.
\v
{\bf 5.} It now remains to prove (ii).
Without loss of generality, we can assume that the intervals $[a_i, b_i]$ are
labelled according to decreasing length, so that
\bel{decl}
b_1-a_1~\geq~b_2-a_2~\geq~\cdots\eeq
 Let $0<\tau<\tau'$ and $\ve>0$ be given.   
Choose 
$\nu>1$ large enough so that (\ref{lsm}) holds.
We now express the open set
$$\,]\tau, \tau'[\,\setminus \left(\bigcup_{1\leq i\leq \nu} [a_i, b_i]
\right)~=~\bigcup_{k=1}^m
\,]\tau_k, \tau'_k[\,$$
as the union of finitely many disjoint open intervals.

Next,
we choose an integer $\nu'>\nu$ such that
\bel{nup}\sum_{i>\nu'} m_1(\Gamma_i)~<~\ve'~\doteq~{\ve\over m}\,,\eeq
and define the times
$$t_k~\doteq~\tau_k+\ve',\qquad t_k'~\doteq~\tau_k'-\ve'.$$
Finally, for $k=1,\ldots,m$, we define the sets of integers
$$I_k~\doteq~\Big\{ i\,;~~\nu+1\,\leq\,i\,\leq\,\nu'\,,~~[a_i, b_i] \cap \,]\tau_k, \tau'_k[~
\not= ~\emptyset\Big\}.$$
Notice that, by (\ref{decl}), these sets are mutually disjoint.
\v
{\bf 6.} Toward a proof of (\ref{lsb}) we will show that, 
for every $k=1,\ldots,m$, one has
\bel{lb3}B\bigl(R^\Gamma(t_k), r_k\bigr) ~\subseteq~\ov{R^\Gamma(\tau_k')},
\qquad\hbox{with}\qquad r_k~=~(\tau_k'-\tau_k )
-3 \ve' - \sum_{i\in I_k} m_1(\Gamma_i).
\eeq
Notice that (\ref{lb3}) implies
\bel{lb4}B\bigl(R^\Gamma(\tau), r\bigr) ~\subseteq~\ov{R^\Gamma(\tau')},
\eeq
with
\bel{r}
\bega{rl}r&=\ds \sum_{k=1}^m r_k~=~\sum_{k=1}^m (\tau'_k-\tau_k-\ve') - 3m \ve' - \sum_{k=1}^m \sum_{i\in I_k} m_1(\Gamma_i)\\[4mm]
&\ds\geq~m_1\left( [\tau, \tau'] \setminus \bigcup_{i=1}^\nu [a_i, b_i] \right) -3m \ve' -
\sum_{\nu<i\leq \nu'} m_1(\Gamma_i)\\[4mm]
&\geq \ds~m_1\left( [\tau, \tau'] \setminus \bigcup_{i=1}^{+\infty} [a_i, b_i] \right) -\ve -3m \ve' -\ve\\[4mm]
&= \ds~m_1\left( [\tau, \tau'] \setminus \bigcup_{i=1}^{+\infty} [a_i, b_i] \right) -5\ve.
\enda\eeq
Since here $\ve>0$ can be taken arbitrarily small, this yields (\ref{lsb}).
\v
{\bf 7.} It thus remains to prove (\ref{lb3}), for each $k\in\{1,\ldots,m\}$.

Consider any point $y\in B(R^\Gamma(t_k, r_k))$ with $y\notin \ov{R^\Gamma(t_k)}$, and 
choose a point $x_0\in\ov{R^\Gamma(t_k)}$ 
which minimizes the distance from $y$.   For any given $\rho>0$, we can then choose
a point $y_0\in R^\Gamma(t_k)$ with $|y_0-x_0|<\rho$.
Notice that we can also assume 
\bel{nden}
\lim_{h\to 0+} {1\over h^2} m_1\bigl(\Gamma\cap B(y_0,h)\bigr)~=~0,\eeq
because this property holds at a.e.~point $x\in \R^2$, w.r.t.~Lebesgue measure.

Call $\gamma$ the segment with endpoints $y_0,y$, and let 
$\xi:[0,\ell]\mapsto \R^2$ be  an arc-length parameterization of this
segment, oriented from $y_0$ to $y$.   Notice that this implies $\ell<r_k$.

In order to use Lemma~\ref{l:22}, we claim that,
among all the connected components $\Gamma_i$, $1\leq i\leq\nu'$
the only ones 
that can have a non-empty intersection with $\gamma$
are the components $\Gamma_i$, with $i\in I_k$.   
Indeed, consider the set of indices
\bel{I-} I^-~\doteq~\bigl\{ i\leq \nu\,;~~
[a_i, b_i]\subseteq [0, \tau_k]\bigr\}~\cup~I_1\cup~\cdots~\cup~I_{k-1}
\eeq
For every $i\in I^-$ we have $b_i\leq \tau_k < t_k$.
Hence
$\ov {R^\Gamma(t_k)}$ contains a neighborhood of $\Gamma_i$. 
Therefore, since $y$ lies outside $\ov {R^\Gamma(t_k)}$, a segment of minimum length joining $y$ 
with a point $x_0\in \ov {R^\Gamma(t_k)}$ cannot intersect $\Gamma_i$.
The same holds if we choose $y_0$ sufficiently close to $x_0$.

Summarizing the previous discussion, given $y\in  B\bigl(R^\Gamma(t_k), r_k\bigr)
\setminus \ov{R^\Gamma(t_k)}$, we can find $y_0\in R^\Gamma(t_k)
$
and a radius $\rho>0$ small enough
such that
\begi
\item[(i)] ~$\ell\,\doteq\,|y_0-y|\,<\, r_k$.
\item[(ii)]  The segment $\gamma$ 
with endpoints  $y_0,y$ does not intersect any of the 
compact connected components $\Gamma_i$ with $i\in I^-$.
\item[(iii)] The 
circumference $\Sigma$ centered at $y_0$ with radius $\rho$ satisfies
\bel{Si}\Sigma~\doteq~\{ x\in\R^2\,;~~|x-y_0| = \rho\}~\subset~R^\Gamma(t_k)\setminus \Gamma.\eeq
\endi
Notice that the (\ref{Si}) is made possible thanks to (\ref{nden}).

Next, consider the set of indices
$$I^+~\doteq~
\bigl\{ i\leq \nu\,;~~a_i\geq \tau'_k\bigr\}\, \cup\, I_{k+1}\,\cup \,\cdots
\,\cup \,I_m\,.$$
Arguing by contradiction,
we show that none of the components $\Gamma_i$ with $i\in I^+$ can 
intersect the segment $\gamma$.   Indeed, if the intersection is nonempty,
define
$$\bar s~\doteq~\min\bigl\{ s\in [0,\ell],;~~\xi(s)\in \Gamma_i\quad
\hbox{for some}~~i\in I^+\bigr\}.$$
For every $s<\bar s$, an application of Lemma~\ref{l:41} would imply the existence of a sequence of paths
$\xi_j:[0, \ell_j]\mapsto \R^2\setminus\Gamma$ 
such that 
 $$\xi_j(0)\to y_0,\qquad \xi_j(\ell_j)\to \xi(s),\qquad\quad\hbox{as}~~j\to\infty,$$
and whose length satisfies
\bel{lj}
\limsup_{j\to\infty} ~\ell_j~\leq~s + \sum_{i\in I_k} m_1(\Gamma_i) + \sum_{i>\nu'}
m_1(\Gamma_i)~<~r_k +  \sum_{i\in I_k} m_1(\Gamma_i) +\ve'\,.\eeq
We now observe that, for all $j$ large enough, the path $\xi_j(\cdot)$ crosses the circumference
$\Sigma$ at some point $\xi_j(s_j)$.   By taking the restriction of $\xi_j$ to the subinterval $[s_j,\ell_j]$, we obtain a sequence of paths $\Tilde\xi_j$, 
of length
$\leq \ell_j$, where the initial point lies on $\Sigma\subset R^\Gamma(t_k)$
and the terminal points converge to $y$.   This implies
$$T^\Gamma(\xi(s))~\leq~t_k + \liminf_{j\to\infty} \ell_j~\leq~t_k+r_k + \sum_{i\in I_k} m_1(\Gamma_i) +\ve'. $$
Since $s$ can be taken arbitrarily close to $\bar s$, recalling (\ref{lb3})
we conclude
 that the point $\xi(\bar s)\in \Gamma_{i^*}$ lies inside the set
 $\ov{R^\Gamma(T)}$, with 
$$\bega{rl} T&=\ds~t_k+r_k + \sum_{i\in I_k} m_1(\Gamma_i) +\ve'
\\[4mm]
&\ds=~(\tau_k+\ve') + \left[(\tau_k'-\ve'-\tau_k) -3\ve'-\sum_{i\in I_k}m_1(\Gamma_i)
\right]+\sum_{i\in I_k} m_1(\Gamma_i) +\ve'\\[4mm]
& <~\tau_k'\,.
\enda$$
Since $i^*\in I^+$, by definition this implies $a_{i^*}\geq \tau_k'$, 
reaching a contradiction. 
\v
{\bf 8.} In view of the previous step, we can now apply Lemma~\ref{l:41}
to each segment with endpoints $y_0$, $\gamma(s)$, for $0<s<\ell=|y-y_0|$.
This yields a sequence of paths $\Tilde \xi_j:[0,\ell_j]\mapsto\R^2$, joining a point 
$x_j\in \Sigma\subset R^\Gamma(t_k)$ with a point $y_j$ which becomes arbitrarily
close to 
$\xi(s)$ as $j\to \infty$.    
All these paths $\Tilde\xi_j$ do not cross $\Gamma$.   Their lengths $\ell_j$ satisfy the uniform bound
$$\ell_j~\leq~s+ \sum_{i\in I_k^+} m_1(\Gamma_i) + \ve'~\leq~r_k +\sum_{i\in I_k^+} m_1(\Gamma_i) + \ve'~\leq~(\tau'_k-t_k) -\ve'.$$
For every $0\le  s<\ell$, this implies
$$\gamma(s)~\in~\ov{ R^\Gamma(\tau_k' -\ve')}~
\subseteq~\ov{ R^\Gamma(\tau_k')}.$$
Letting $s\to \ell$, we obtain $\gamma(s)\to y$, 
and   hence $y\in \ov{ R^\Gamma(\tau_k')}$.
This establishes the inclusion (\ref{lb3}) for every $k=1,\ldots,m$, thus
completing the proof.
\endproof

\section{Properties of the minimum time function with obstacles}
\label{s:3}
\setcounter{equation}{0}
Assume that the initial set $R_0$ where the fire is burning at $t=0$ has finite perimeter.
Consider a barrier $\Gamma=\cup_i\Gamma_i$, written as  
the union of its connected components.  For every fixed time $\ov T>0$, the truncated function 
$$x~\mapsto ~\min\bigl\{\ov T, T^\Gamma(x)\bigr\}$$ has bounded variation.
Indeed, as shown in \cite{DLR}, it is an SBV function.
By the co-area formula it thus follows
\bel{coa}\int_0^{\ov T} m_1\Big( \partial \ov{R^\Gamma(t)}\Big)\, dt~<~\infty.\eeq
As a consequence, for a.e.~time $t\in [0, \ov T]$, the boundary 
$\partial \ov{R^\Gamma(t)}$
is a curve with finite length. 
\v

We now consider a sequence of barriers $\Gamma^{(n)}$, $n\geq 1$,
converging to 
a barrier $\Gamma$, and study the behavior of the corresponding 
minimum time functions $T^{\Gamma^{(n)}}$.   Two cases will be studied.
The first lemma deals with the case where each barrier has a finite number of connected components. The second lemma is concerned with barriers having
countably many components. For the definition and properties of the Hausdorff distance between compact sets we refer to \cite{AC, BPi}.

\begin{lemma} \label{l:26} Let a bounded open set $R_0\subset\R^2$
be given.
Consider a barrier $\Gamma=\cup_{i=1}^N \Gamma_i$ which is the union of finitely many
compact, simply connected, rectifiable components.
Let $\Gamma^{(n)}= \cup_{i=1}^N \Gamma_i^{(n)}$, with $n\geq 1$,
be an approximating sequence of barriers.
Assume the convergence w.r.t.~the Hausdorff distance:
\bel{Hconv}
\lim_{n\to\infty} ~d_H(\Gamma_i^{(n)}, \Gamma_i)~=~0\qquad\hbox{for each}~ i=1,\ldots,N.
\eeq
\begi
\item[(i)] For every $x\in\R^2\setminus \Gamma$, one has
\bel{ltn}
T^\Gamma(x)~=~\lim_{n\to\infty}~ T^{\Gamma^{(n)}}(x).\eeq
\item[(ii)] For each $n\geq 1$, let $ \xi_n:[0, \tau_n]\mapsto \R^2$ be an optimal trajectory reaching a point $x\in \R^2\setminus\Gamma$ in minimum
time without crossing $\Gamma^{(n)}$. If  $\tau_n\to \ov \tau$ and $\xi_n(\cdot)\to \xi(\cdot)$ uniformly on every compact subset of $[0, \ov \tau[\,$,  then $\xi(\cdot)$ 
is an optimal trajectory reaching $x$ in minimum time without crossing 
$\Gamma$.
\endi
\end{lemma}

{\bf Proof.}   {\bf 1.} Toward a proof of 
 (i), consider a minimizing sequence of 1-Lipschitz paths $\xi_\nu:[0, \tau_\nu]\mapsto \R^2$  , satisfying
$$ \xi_\nu(0)\in R_0,\qquad \xi_\nu(t)\notin \Gamma \quad\forall t\in[0,\tau_\nu],$$
$$|\xi(\tau_\nu)-x|<{1\over \nu}\,,\qquad\qquad \tau_\nu~\to~T^{\Gamma}(x) \quad
\hbox{as}\quad \nu\to\infty.$$
Since $\Gamma$ is compact and by assumption $x\notin \Gamma$, we conclude that, for all $\nu$ large enough, the segment joining $\xi(\tau_\nu)$ with $x$
will not touch $\Gamma$.    By adding this segment to the path $\xi_\nu$
we obtain another sequence of 1-Lipschitz paths $\tilde\xi_\nu:[0, \Tilde T_\nu]\mapsto\R^2\setminus \Gamma$, with 
\bel{TTn}\tilde\xi_\nu(\Tilde T_\nu) = x \qquad \forall \nu\,,\qquad\qquad \lim_{\nu\to\infty}
\Tilde T_\nu~=~T^{\Gamma}(x) .\eeq

By the compactness of the barriers $\Gamma^{(n)}$, for each $\nu\geq 1$ there exist an integer $n_\nu$ large enough such that   $\tilde 
\xi_\nu(t)\notin \Gamma^{(n)}$ for all 
$n \geq n_\nu$ and all 
$t\in [0, \Tilde T_\nu]$. This immediately implies 
$$\Tilde T_\nu~\geq~\limsup_{n\to\infty}~ T^{\Gamma^{(n)}}(x).$$
Together with (\ref{TTn}), this yields
\bel{lt4}
T^\Gamma(x)~\geq ~\limsup_{n\to\infty}~ T^{\Gamma^{(n)}}(x).\eeq

\v
{\bf 2.} To prove (ii) we observe that, by (\ref{lt4}), 
$\ov \tau\leq T^\Gamma(x)$.    It thus only remains to prove that 
the limit path $\xi(\cdot)$ is admissible.  Since $x\notin \Gamma$,
there exists $\rho>0$ such that $B(x,\rho)\cap\Gamma= \emptyset$.
In the following, to simplify notation, we still denote by $\xi$ the set of points
$\{ \xi(t)\,;~t\in [0,\ov\tau]\}\subset\R^2$. For 
a given $0<r<\!<\rho$
sufficiently small, consider the neighborhood $B(\xi, r)$, and let
$\Sigma$ be the boundary of the unbounded connected component
of $\R^2\setminus B(\xi,r)$.    This is a simple closed curve, which we
can parameterize by arc-length,  oriented
counterclockwise.   

\begin{figure}[htbp]
   \centering
 \includegraphics[width=0.7\textwidth]{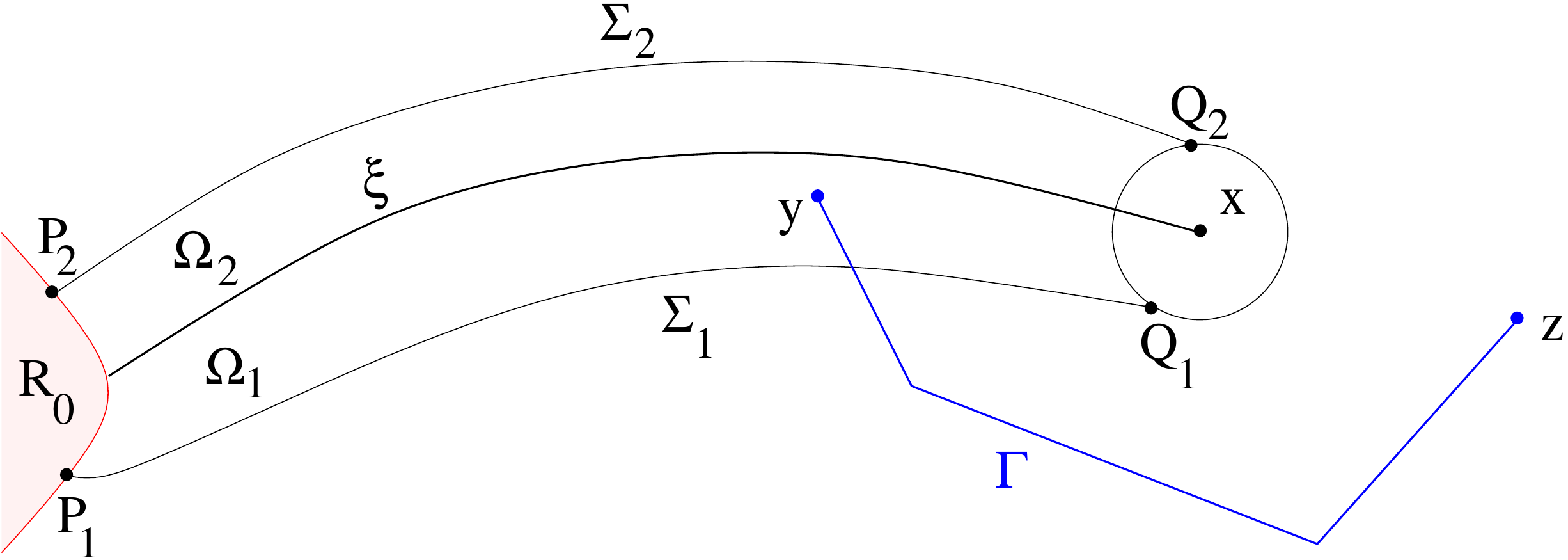}
   \caption{\small Showing that the limit path $\xi$ is admissible.}\label{f:fc197}
\end{figure}

As shown in Fig.~\ref{f:fc197}, 
within $\Sigma$ we distinguish an arc $\Sigma_1$ connecting a point in 
$P_1\in R_0$ with a point $Q_1\in B(x, \rho)$, and an arc
$\Sigma_2$ connecting a point $Q_2\in  B(x, \rho)$ with a point
$P_2\in  R_0$, moving counterclockwise.  
we call $\Hat\Omega$ the open set bounded by 
$\Sigma$. Moreover, we consider the open subset
$$\Omega~\doteq~ \Hat \Omega\setminus (\ov R_0 \cup \ov B(x,\rho))$$
with its open subsets
$\Omega_1\subset\Omega$,  bounded between $\Sigma_1$ and $\xi$,
 and  $\Omega_2\subset\Omega$
 bounded between 
$\xi$ and $ \Sigma_2$.
\v
{\bf 3.}
In the following, for simplicity we consider the case where $N=1$, so that 
$\Gamma$ and all the approximating barriers $\Gamma^{(n)}$ contain only one component.
Since all components are compact and have a positive distance from each other, the general case follows by the same arguments.

Fix a point $z\in \Gamma$, outside the region enclosed by $\Sigma$.
We claim that, for every $y\in \Omega_1\cap\Gamma$, 
there exists a path $\gamma^y\subseteq \Gamma$ connecting
$y$ with $z$, and touching $\Sigma_1$ without entering $\Omega_2$.
More precisely, we claim that there is a map  
$\gamma^y:[0, \bar s]\mapsto \Gamma$ such that
$$\gamma^y(0)=y,\quad \gamma^y(\bar s) = z,$$
$$\gamma^y(s^*) ~\in \Sigma_1,\quad\hbox{where}
\quad s^*\doteq\inf~ \Big\{s\in [0, \bar s]\,;~~\gamma^y(s)\notin \Omega\Big\}.$$
To prove this claim, we observe that there exist
sequences $y_n, z_n\in \Gamma^{(n)}$, with 
$y_n\to y$ and $z_n\to z$.  Since $\Gamma^{(n)}$ is connected, 
for each $n\geq 1$ there is a path $\gamma_n$ joining $y_n$ with 
$z_n$, and remaining inside $\Gamma^{(n)}$.  
By possibly selecting a subsequence and relabeling, we obtain a limit path
$\gamma:[0, \bar s]\mapsto \Gamma$, 
joining $y$ with $z$.   If $\gamma(s)\in \Omega_2$ for some $0<s<s^*$,
then $\gamma$ crosses the path $\xi$.   By uniform convergence
$\gamma_n\to \gamma$ and $\xi_n\to \xi$, this would imply that 
every $\xi_n$, with $n$ suitably large,
crosses $\Gamma^{(n)}$, a contradiction.

We observe that, after reaching the boundary $\Sigma_1$, for $s\in [s^*, \bar s]$ the path $\gamma^y$ can re-enter inside $\Omega$.  However, it cannot
cross $\xi$.  Namely, if it enters through $\Sigma_1$, it must eventually
leave through $\Sigma_1$.   If it enters through $\Sigma_2$, it must 
leave through $\Sigma_2$. Otherwise, being a limit of paths $\gamma_n$ contained in 
the approximating barriers $\Gamma^{(n)}$, these paths would cross
the corresponding paths $\xi_n$.   

\v

{\bf 4.} Within the compact curve $\xi$, for $k=1,2$ we  define the subset
\bel{Vk}
V_k~\doteq~\Big\{ y\in \xi\cap\Gamma\,;~~\hbox{there
is a path inside $\Gamma$ joining $y$ with $z$, exiting through $\Sigma_k$}
\Big\}.\eeq
Since $\Gamma$ is rectifiable and compact, it follows that $V_1, V_2$ are both compact.
We claim that they are disjoint.
Indeed, assume that $y=\xi(s) \in V_1\cap V_2$.
Then, inside $\Omega\cap\Gamma$, we can find a path joining $y$ with a point $y_1\in \Omega_1$, and another path joining $y$ with a point 
$y_2\in \Omega_2$.   In turn, by the previous step, there is a path
$\gamma^{y_1}$ joining $y_1$ to $z$, and a path $\gamma^{y_2}$
joining $y_2$ with $z$.    The union of these paths is a multiply
connected rectifiable subset of $\Gamma$.   This yields a contradiction.
\v
{\bf 5.} By the previous step, we can cover the disjoint compact sets
$V_1, V_2\subset [0, \ov\tau]$ with finitely 
many disjoint intervals, say
$[a_j, b_j]$ and $[c_j, d_j]$, so that
$$V_1~\subseteq~\bigcup_{j=1}^m [a_j, b_j],\qquad\qquad 
V_2~\subseteq~\bigcup_{j=1}^m [c_j, d_j].$$
Define the corresponding  portions of curve
$$\gamma_1^{(j)}~\doteq~\bigl\{\gamma(s)\,;~~s\in [a_j, b_j]\bigr\},
\qquad \gamma_2^{(j)}~\doteq~\bigl\{\gamma(s)\,;~~s\in [c_j, d_j]
\bigr\}.$$
We claim that, by choosing  a radius $\delta>0$ small enough, 
for every $j=1,\ldots, m$ one has
\bel{noni}
\Gamma\cap B\bigl(\gamma_1^{(j)}, \delta\bigr)\cap \Omega_2~=~\emptyset,
\qquad\qquad
\Gamma\cap B\bigl(\gamma_2^{(j)}, \delta\bigr)\cap \Omega_1~=~\emptyset.
\eeq
Indeed, if no such radius $\delta>0$ exists, we could find a point 
$y\in V_1$ and a sequence of points $y_n\to y$ with 
$y_n\in \Gamma\cap \Omega_2$ for all $n\geq 1$.
By step 3, for each $y_n$ there exists a  path joining $y_n$ to 
$z$, remaining inside $\Gamma\setminus \Omega_1$.
By taking a limit, we obtain a path joining $y$ with $z$, remaining inside 
$\Gamma\setminus \Omega_1$.   This would yield $y\in V_2$, 
reaching a contradiction because in step 2 we proved that $V_1\cap V_2=\emptyset$.
\v
\begin{figure}[htbp]
   \centering
 \includegraphics[width=0.7\textwidth]{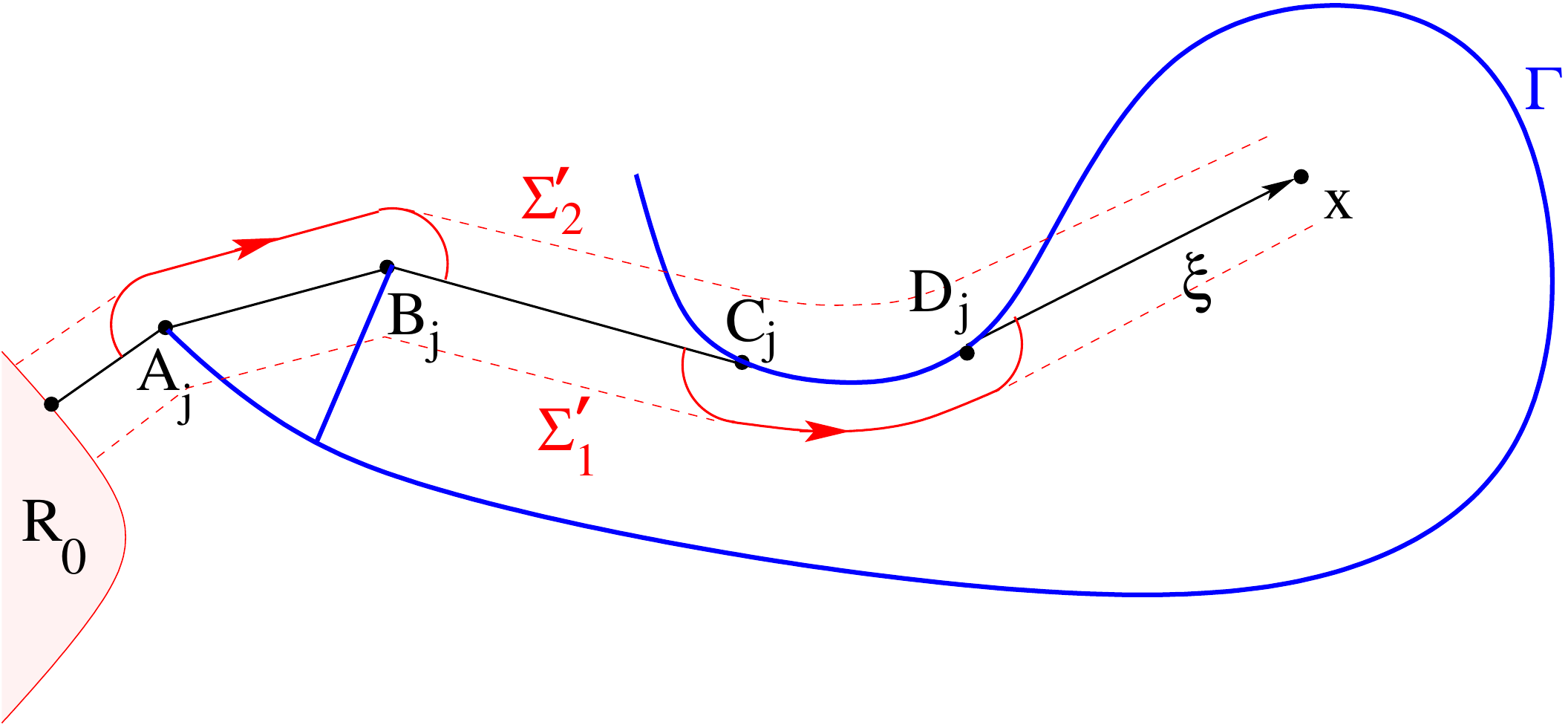}
   \caption{\small By modifying the trajectory $\xi(\cdot)$ along the arcs where it touches the barrier $\Gamma$, one can show that $\xi$ is admissible.}\label{f:fc198}
\end{figure}

{\bf 6.} We now describe how to make a small modification of the path $\xi$, so that it does not touch the barrier $\Gamma$.  Fix $\ve>0$ and consider the finitely many circumferences with radius $\ve$, centered at the points 
$$A_j= \xi(a_j),\quad B_j = \xi(b_j), \qquad 
C_j= \xi(c_j), \qquad D_j = \xi(d_j).$$
In addition, for $0<\ve'<\!<\ve$, call $\Sigma'$ the simple closed curve
obtained by taking the boundary of the unbounded connected component
of $\R^2\setminus B(\xi,\ve')$.  As in step  2, we distinguish a lower and an upper portion of this boundary, which we call $\Sigma'_1, \Sigma'_2$, respectively.

 As shown in Fig.~\ref{f:fc198}, for each   $j=1,\ldots,m$, 
 the portion of the path 
 $\{\xi(s)\,;~s\in [a_j, b_j]\}$ between 
 $A_j$ and $B_j$, is replaced by two arcs of circumferences
 centered at $A_j, B_j$ together with a portion of the curve $\Sigma_2'$.
 Similarly, the portion of the path 
 $\{\xi(s)\,;~s\in [c_j, d_j]\}$ between 
 $C_j$ and $D_j$, is replaced by two arcs of circumferences
 centered at $C_j, D_j$ together with a portion of the curve $\Sigma_1'$.
By the previous analysis, for all $\ve,\ve'>0$ sufficiently small, this new curve
does not intersect $\Gamma$.   Moreover, letting $\ve,\ve'\to 0$, 
we recover the original 
path $\xi$ in the limit.   This shows that $\xi(\cdot)$ is admissible, proving (ii).
\v
{\bf 8.}   By (ii) it now follows
$$T^\Gamma(x)~\leq~\liminf_{n\to\infty} T^{\Gamma^{(n)}}(x).$$
Together with (\ref{lt4}), this yields (\ref{ltn}), completing the proof.
\endproof

\begin{figure}[htbp]
   \centering
 \includegraphics[width=0.7\textwidth]{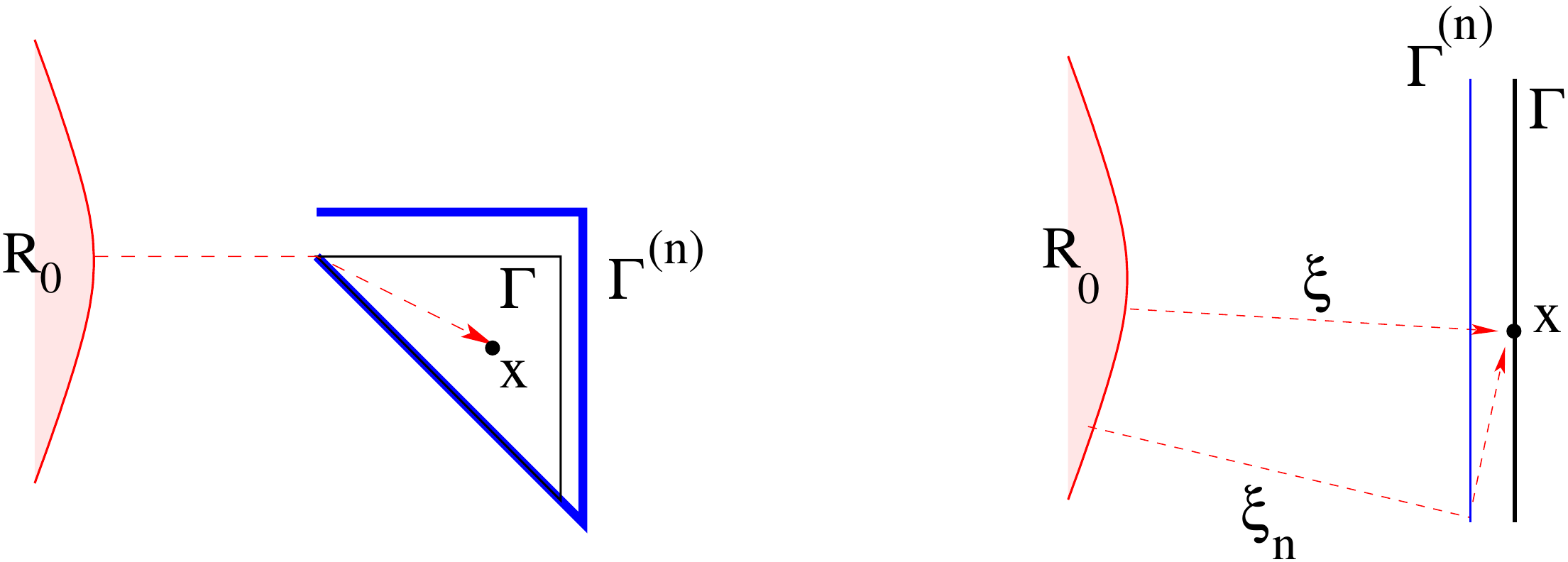}
   \caption{\small Left: an example showing that  (\ref{ltn}) can fail, if the 
   barrier $\Gamma$ is not simply connected.
   For each $n\geq 1$, the barrier $\Gamma^{(n)}$ is the union of three segments, and $\R^2\setminus\Gamma^{(n)}$ is connected.   However, the limit barrier $\Gamma$ is the boundary of a triangle, which is not simply connected.  None of the points $x$ in the interior of this triangle can 
   be reached from $R_0$, without crossing $\Gamma$. 
   Right: an example showing that, if $x\in\Gamma$, the inequality (\ref{ltn}) can fail. Indeed, here $T^\Gamma(x)< \lim_{n\to\infty} T^{\Gamma^{(n)}} (x)$. }\label{f:fc185}
\end{figure}

\begin{remark} {\rm  In the above lemma, the assumptions that each $\Gamma_i$ is simply connected and that $x\notin \Gamma$ play an essential role.  In Figure~\ref{f:fc185} shows two cases where these assumptions
are not satisfied, and the conclusions fail.
}
\end{remark}

\begin{remark} {\rm In (\ref{Vk}),
one can think of $V_1$ is the set of points where the barrier $\Gamma$
touches the optimal trajectory $\xi$ on the right, while $V_2$ is the set of points where $\Gamma$
touches $\xi$ on the left.
Calling $\dot\xi(t)=(\cos\theta, \sin\theta)\in \R^2$ the tangent vector,
by construction the map $t\mapsto \theta(t)$ is non-increasing along each 
interval $[a_j, b_j]$,  non-decreasing along each interval $[c_j, d_j]$,
and constant everywhere else.}
\end{remark}

\v

\begin{lemma} \label{l:27} Let a bounded open set $R_0\subset\R^2$
be given.
Consider a barrier $\Gamma=\cup_{i=1}^\infty \Gamma_i$ and assume that 
$\R^2\setminus\Gamma$ is connected.
For each $\nu\geq 1$, consider the finite union $\Gamma_\nu= \cup_{i=1}^\nu \Gamma_i$.   Call $T^\Gamma$, $T^{\Gamma_\nu}$ the corresponding
minimum time functions.
\begi
\item[(i)] For every $x\in\R^2$ one has
\bel{ltnu}
T^\Gamma(x)~=~\lim_{\nu\to\infty} T^{\Gamma_\nu}(x).\eeq
\item[(ii)] For each $\nu\geq 1$, let $ \xi_\nu:[0, \tau_\nu]\mapsto \R^2$ be an optimal trajectory reaching $x$ in minimum
time without crossing $\Gamma_\nu$. If  $\tau_\nu\to \ov \tau$ 
and $\xi_\nu(\cdot)\to \xi(\cdot)$ uniformly on every compact subset of 
$[0, \ov \tau[\,$,  then $\xi(\cdot)$ 
is an optimal trajectory reaching $x$ in minimum time without crossing~$\Gamma$.
\endi
\end{lemma}

{\bf Proof.}   
{\bf 1.}
To prove (\ref{ltnu}), fix $x\in \R^2$ and, for every $\nu\geq 1$, call $\tau_\nu\doteq T^{\Gamma_\nu}(x)$.   Denote by $\zeta_\nu:[0,\tau_\nu]\to \R^2$ an optimal trajectory reaching the point $x$ without crossing $\Gamma_\nu$. 
According to  Definition~\ref{d:adm}, there exists a second path
$\xi_\nu:[0,\tau_\nu]\to \R^2\setminus \Gamma_\nu$ such that 
$$\xi_\nu(0)\in R_0\,,\qquad |\xi(t)-\zeta(t)\bigr|~<~{1\over\nu}\qquad
\forall t\in [0, \tau_\nu].$$
Applying 
Lemma~\ref{l:22}, we obtain a further path $\tilde\xi_\nu:[0,\Tilde\tau_\nu]\to\R^2$, also parameterized by arc length, such that 
$$\tilde \xi(0)\in R_0\,,
\qquad \qquad \bigl|\tilde\xi_\nu(\Tilde \tau_\nu)-x\bigr|~\leq~{2\over \nu},$$
$$\tilde\xi_\nu(t)\notin\Gamma\qquad\forall t\in [0,\Tilde \tau_\nu],$$
and with length
$$\Tilde \tau_\nu~<~\tau_\nu+ \sum_{i>\nu} m_1(\Gamma_i)\,.$$
Therefore (\ref{ltnu}) follows from
$$ \limsup_{\nu\to \infty}\tau_\nu~\leq~ T^\Gamma(x) ~\leq ~\liminf_{\nu\to \infty}
\Tilde\tau_\nu
~=~\liminf_{\nu\to\infty} \left[\tau_\nu+ \sum_{i>\nu} m_1(\Gamma_i)\right]~=~\liminf_{\nu\to\infty} \tau_\nu\,. $$

\v
{\bf 2.}
To prove part (ii), as usual we assume  that all the optimal trajectories $\xi_\nu$ are parameterized by arc length. By the previous step one has
$$ T^{\Gamma}(x) ~=~ \lim_{\nu\to\infty} T^{\Gamma_\nu}(x)~ =~
 \lim_{\nu\to\infty} \tau_n ~=~\ov \tau. $$
To achieve the proof it thus suffices to check that the limit trajectory
$\xi(\cdot)$ is admissible. 

Toward this goal, the key tool is again provided by Lemma \ref{l:22}. 
For each $\nu\geq 1$, using the lemma we obtain 
a path $\tilde{\xi}_\nu:[0,\Tilde{\tau}_\nu]\to \R^2$ such that 
$$\tilde{\xi}_\nu(0) \in R_0\,,
\qquad\qquad \tilde{\xi}_\nu(t)\notin\Gamma\qquad\forall t\in [0,\Tilde{\tau}_\nu],
$$
$$
\bigl|\tilde{\xi}_\nu(\Tilde\tau_\nu)-\xi_\nu(\tau_\nu)\bigr|~\leq~{1\over \nu},
$$
and with length
$$\Tilde{\tau}_\nu~\leq~\tau_\nu+ \sum_{i>\nu} m_1(\Gamma_i)\,.$$

Recalling (\ref{xijj}) in the proof of Lemma~\ref{l:22}, w.l.o.g.~we can assume
that 
\bel{unif}
\bigl|\tilde \xi_\nu (t) -\xi_\nu(t)\bigr|~\leq~\left(1+{1\over \nu}\right)\sum_{i>\nu} m_1(\Gamma_i),\eeq
for all $t\geq 0$.  It is understood that here  $\xi_\nu$ and $\tilde\xi_\nu$ 
are extended as   constant functions,
for $t\geq \tau_\nu$ and $t\geq \Tilde\tau_\nu$, respectively.

By (\ref{unif}), as $\nu\to \infty$ the sequence of paths $\tilde{\xi}_\nu(\cdot)$ converges uniformly to $\xi(\cdot)$.  Hence
$\xi$ is admissible.

\section{A regularity property of optimal trajectories}
\label{s:4}
\setcounter{equation}{0}

Aim of this section is to study a property of
the optimal trajectories for the fire, in the presence of barriers. We begin with a few observations.
\begi
\item If no barriers are present, all optimal trajectories are straight lines, and the 
minimum time function is trivially $T(x)= d(x, R_0)$.

\begin{figure}[htbp]
   \centering
 \includegraphics[width=0.9\textwidth]{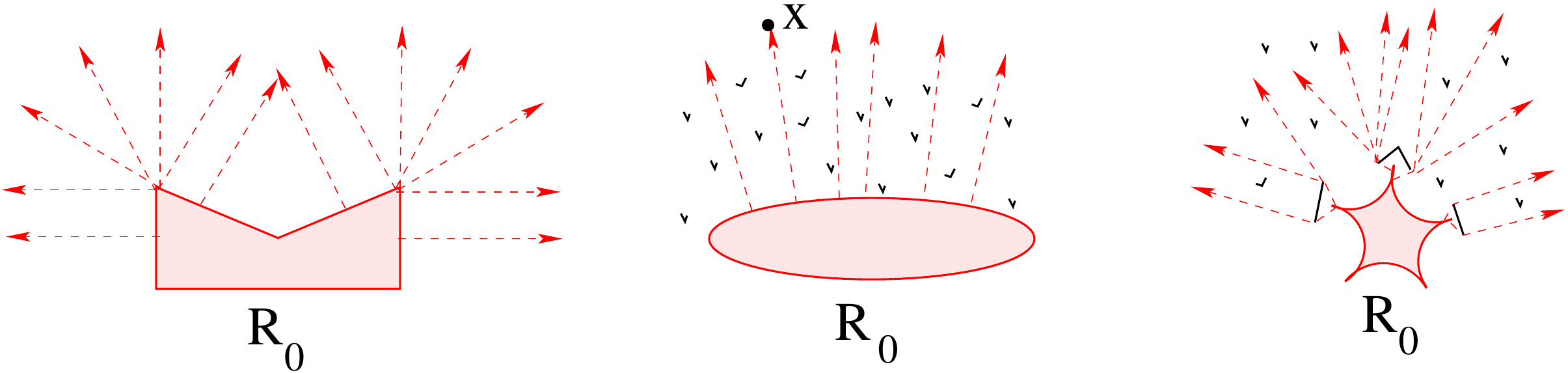}
   \caption{\small Left: if no barriers are present, all optimal trajectories are straight lines.
   Center: if the initial set $R_0$ has smooth boundary and the 
   total length of all barriers is small,   then most of the optimal trajectories 
   do not hit the barrier, and hence are still straight lines.
    Right: for a general set $R_0$ whose boundary contains cusps, even if the total length of the barriers is small,
    most of the optimal trajectories may touch one of the barriers.}
   \label{f:fc175}
\end{figure}

\item Next, assume that $R_0$ has a $\C^2$ boundary. For each point $x\notin R_0$,
consider the shortest segment connecting $x$ with a point $y\in \ov{R_0}$. 
If the total length
of all barriers $m_1(\Gamma)= \sum_i m_1(\Gamma_i)<\ve$ is sufficiently small,
then most of these segments will not cross $\Gamma$.   Hence, as shown in Fig.~\ref{f:fc175}, center,  they will yield 
optimal trajectories for the fire also when barriers are present.
Notice that this remains true even if the set $\Gamma$ of all barriers is dense in $\R^2$.

\item For a general set $R_0$, however, even if the total length of all barriers is very small,
it can happen that most of the optimal trajectories touch one of the barriers.
As shown in Fig.~\ref{f:fc175}, this is the case when $R_0$ has  cusps, and some of 
the barriers are placed very close to these cusps.
\endi

Since in general it is not true that most optimal trajectories are straight lines, in this section we prove a somewhat weaker property. Namely: most optimal trajectories contain 
long straight segments.  This property will play a key role in the proof of 
Theorem~\ref{t:1}.

Consider again the optimization problem for the fire, in the presence of a 
barrier $\Gamma$. Given an initial set 
$R_0$, let  
\bel{R1def}R_1~\doteq ~\{x\in\R^2\,;~d(x,R_0)<1\}\eeq 
be the neighborhood of radius 1
around $R_0$.
For any $x\in \R^2$, call  $T^\Gamma(x)$ the minimum time needed to reach $x$ from $R_0$
without crossing  $\Gamma$.    
Moreover, given an admissible trajectory $t\mapsto \xi^x(t)$ reaching $x$ in minimum time,
we denote by $\rho(x)$  the length of the last portion of this trajectory which is
a straight line.
More precisely, 
\bel{rhodef}\bega{rl} \rho(x)&\doteq~\sup~\Big\{\tau\geq 0\,;~~
\hbox{there exists a trajectory $t\mapsto \xi^x(t)$ reaching $x$ in minimum time 
} \\[4mm]& \qquad\quad \hbox{without crossing $\Gamma$, and  the velocity $\dot \xi^x $ is constant
on }~[T^\Gamma(x)-\tau, ~T^\Gamma(x)]\Big\}.\enda\eeq
\v

\begin{lemma}
\label{l:5}
Let $R_0$ be a bounded, open set, and call $R_1$ the set in (\ref{R1def}).
Then, for any barrier $\Gamma$ one has
%
\bel{rb0}
\int_{R_1} \rho(x)\, dx ~\geq~\int_{R_1} d(x,R_0)\, dx - {\Hat T^2+\Hat T\over 2}\cdot m_1(\Gamma),\eeq
where
$$\Hat T~\doteq~\sup_{x\in R_1} T^\Gamma(x).$$
\end{lemma}

\begin{remark} {\rm In the case where no barriers are present,
one has $\rho(x)= d(x, R_0)$ and the bound (\ref{rb0}) is obvious.
We observe that  a lower bound on the left hand side of (\ref{rb0}) cannot be achieved by the trivial estimate
\bel{gg}\rho(x)~\geq~\inf_{y\in\Gamma} |y-x|,\eeq
because the set $\Gamma=\cup_i\Gamma_i$ can be everywhere dense. In this case  the right hand side of (\ref{gg}) is identically zero.}
\end{remark}

\begin{remark} {\rm Assuming that $\R^2\setminus \Gamma$ 
is connected, so that all barriers are 
only delaying the fire, by Lemma~\ref{l:22} it follows that 
\bel{HTB}\Hat T~\leq~1+m_1(\Gamma).\eeq
}
\end{remark}

\v
\subsection{Polygonal barriers.}
We shall give a proof of   Lemma~\ref{l:5} first in a special case where
explicit computations can be performed.  The general case will then be handled by an
approximation argument. 
In this section, we tassume
\begi
\item[{\bf (A2)}] {\it 
The initial set  $R_0$ is the union of finitely many open discs, while the barrier $\Gamma$ is the union of finitely many (not necessarily disjoint) closed segments. } 
\endi
Notice that this special setting implies
\begi
\item[(i)] Every optimal trajectory for the fire, 
reaching a point $x\in\R^2$ in minimum time without crossing the barriers,
is a polygonal, say with vertices $P_0, P_1,\ldots P_N$.
Here $P_0\in \ov {R_0}$, $P_N=x$, while $P_i\in\Gamma$ for all $i\in \{1,\ldots, N-1\}$. Indeed,
 each $P_i$ will be an edge of one of the segments forming the barrier $\Gamma$.
\item[(ii)] For every $t>0$, the boundary of the reachable set $\partial \ov{R^\Gamma(t)}$
is the union of finitely many arcs of circumferences.  
\item[(iii)] The set of points which can be reached in minimum time by two
distinct trajectories is the union of finitely many segments, or arcs of hyperbolas.  
\endi
To prove the estimate (\ref{rb0})
we shall study a family
of problems, parameterized by time.
Call
$$\Gamma(t)~=~\Gamma\cap\ov {R^\Gamma(t)}$$
the portion of the walls which are touched by the fire within time $t$.  
We obviously have
$$\Gamma(s)~\subseteq~\Gamma(t)
\qquad\hbox{for} ~~ s<t.$$
For every $t\geq 0$, call $\rho(t,x)$ the function defined at (\ref{rhodef}), but with 
$\Gamma$ replaced by the smaller set $\Gamma(t)$.

\begin{figure}[htbp]
   \centering
 \includegraphics[width=0.9\textwidth]{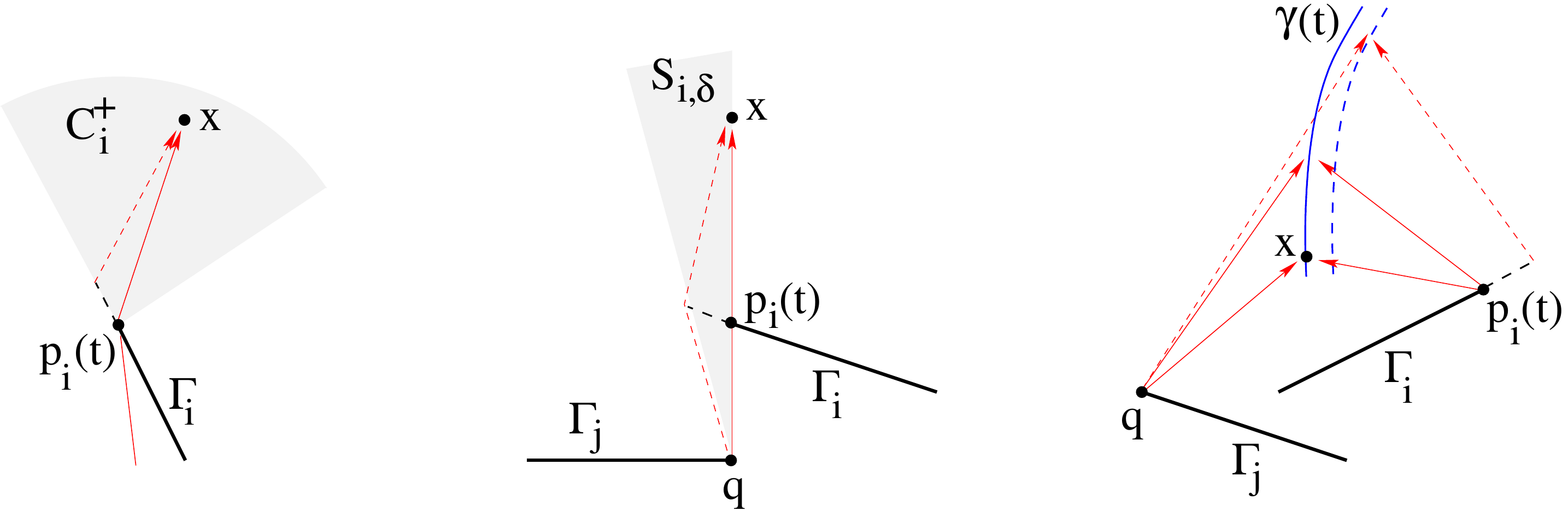}
   \caption{\small Left: as  the segment $\Gamma_i$ becomes longer, the value $\rho(x)$ can decrease at all points $x$ in the shaded region.
   Center: at time increases, the value $\rho(x)$ jumps downward from $|x-q|$ to 
   $|x-p_i(t)|$.  Right: the curve $\gamma(t)$ denotes the set of points $x$ reached in 
   minimum time by two distinct trajectories.   As time increases, this curve changes in time. At points $x$ on this curve, the value of $\rho$ jumps upward from $|x-p_i(t)|$
   to $|x-q|$. 
   }
   \label{f:fc176}
\end{figure}

The estimate (\ref{rb0}) will be achieved by showing that, for a.e.~$t\in [0,1]$,
\bel{c1}
- {d\over dt} \int_{R_1}\rho(t,x)\, dx~\leq~ {\Hat T^2+\Hat T\over 2}\cdot {d\over dt} m_1(\Gamma(t)).\eeq

To fix the ideas, let $\Gamma_i(t)\subseteq\Gamma(t)$ be one of the segments of the barrier, with an endpoint $p_i(t)\in \partial \ov{R^\Gamma(t)}$ moving along the edge
of the advancing fire.    
For a fixed time $t$, referring to the optimization problem with barrier $\Gamma(t)$,
three cases must be considered.
\v
{\bf 1.} Points $x$  reached in minimum time by a trajectory which touches the
point $p_i(t)$.

The set  of all these points, that we shall call $\Omega_i$,  
is contained within a half disc
$C_i$, with center at $p_i(t)$ and radius $\Hat T-t$.
As time increases, 
 for all $x\in \Omega_i$ we have
\bel{t1}
{d\over dt} \rho(t,x)~=~{d\over dt} |p_i(t)-x|~=~ \left\langle \dot p_i(t),\, {p_i(t)-x\over |p_i(t)-x|}\right\rangle~\geq~-\bigl|\dot p_i(t)\bigr|.
 \eeq
Observe that the quantity in (\ref{t1}) can be negative only on the quarter disc
$$C_i^+~\doteq~\bigl\{ x\in C_i\,;~\langle  \dot p_i(t),\, p_i(t)-x\rangle <0\bigr\},
$$
corresponding to the 
shaded region in Fig.~\ref{f:fc176}, left.  
Using (\ref{t1}) we compute
\bel{t3}\bega{l}\ds
{d\over dt}\int_{\Omega_i}\rho(t,x)\, dx~ \geq ~\int_{C_i^+} \left\langle \dot p_i(t),\, {p_i(t)-x\over |p_i(t)-x|}\right\rangle\, dx\\[4mm]
\qquad\ds
\geq~
-\bigl|\dot p_i(t)\bigr|\cdot \int_0^{\Hat T-t} \int_0^{\pi/2} \cos\theta\, d\theta\, r \, dr
~\geq~-{\Hat T^2\over 2}\,\bigl|\dot p_i(t)\bigr|\,.\enda
\eeq
\v
{\bf 2.} Next, we consider points $x$ reached in minimum time by a trajectory whose last 
portion is a segment with endpoints $q$ and $x$, and such that $p_i(t)$ is a point
inside this segment (see Fig.~\ref{f:fc176}, center).

The set of all these points, which we will call $D_i$, is contained on a half line starting at $q$ and passing through $p_i(t)$, so that
\[ |x-p_i(t)|~\leq~|x-q|~<~1\qquad \forall x\in D_i\,.\]
As time increases, 
the value of $\rho$ along $D_i$ jumps downward from $|x-q|$ to $|x-p_i(t)|$. 
To compute the rate of decrease in  the integral $\int\rho(t,x)\, dx$ due to such points, fix $\delta>0$ small and consider the region $D_{i,\delta}$ of all points $x$ such that 
\[\rho(t,x)\,=\,|x-q| \,,\qquad \quad\rho(t+\delta,x)\,=\,|x-p_i(t+\delta)|.\]
By the triangle inequality one obtains
\bel{j2}
\int_{D_{i,\delta}}\bigl(\rho(t,x)-\rho(t+\delta,x)\bigr)\, dx~\leq ~m_2(D_{i,\delta})\cdot 
\bigl|p_i(t+\delta)-q\bigr|.
\eeq

Observe that $D_{i,\delta}$ is contained in a circular sector $S_{i,\delta}$
with radius $\Hat T$  (as the one shaded in Fig.~\ref{f:fc176}, center) whose area can be computed using 
the vector product
\bel{j3}
m_2(S_{i,\delta})~=~\frac{\Hat T}{2} \frac{\Big| 
\bigl(p_i(t+\delta)-p_i(t)\bigr)\times \bigl(p_i(t)-q\bigr)\Big|}
{\bigl|p_i(t) - q\bigr|}+o(\delta).
\eeq
Combining (\ref{j2}) with (\ref{j3}),  we conclude
\bel{j4}\bega{rl} 
 \ds\lim_{\delta\rightarrow 0^+}{1\over\delta}\int_{D_{i,\delta}}\bigl(\rho(t,x)-\rho(t+\delta,x)\bigr)\, dx
&\ds\leq ~\lim_{\delta\rightarrow 0^+}{1\over\delta}\, m_2(S_{i,\delta})\cdot|p_i(t+\delta)-q|~ \leq ~\frac{\Hat T}{2}\, \bigl|\dot p_i(t)\bigr|.
\enda \eeq
\vs

{\bf 3.} Points $x$ on a curve $\gamma(t)$ reached in minimum time by 
two distinct optimal trajectories. 

As shown in Fig.~\ref{f:fc176}, right, we can assume that one of these touches $p_i(t)\in \Gamma_i$, while the other 
touches some other point $q\in\Gamma_j$. These two trajectories have the same length, therefore
\bel{ul}
|x-q|+T^{\Gamma}(q)~=~|x-p_i(t)|+T^{\Gamma}(p_i(t))~=~|x-p_i(t)|+t.
\eeq
Since $T^{\Gamma}(q)\leq t$, this implies $|x-q|\geq |x-p_i(t)|$ for all $x\in\gamma(t)$.

Notice that $\gamma(t)$ is a branch of hyperbola. For a.e.~$y\in \gamma(t)$
we can choose a neighborhood $V$ of $y$ such that, for every $x\in V$, 
one has
\bel{mit}
T^{\Gamma(t)}(x)~=~\min\Big\{ |x-q| + T^\Gamma(q),\quad |x-p_i(t)| + t\Big\}.\eeq
Assume  that, when the barrier is $\Gamma(t)$,
a point  $x\in V$ is reached in minimum time by a trajectory passing through $q$, namely
$$|x-q|+ T^\Gamma(q)~<~|x-p_i(t)|+t.$$
As time increases from $t$ to $t+\delta$ and the point 
 $p_i(t)$ is replaced by $p_i(t+\delta)$, by (\ref{mit}) we have
 \bel{ch}\bega{l}
\bigl|x-p_i(t+\delta)\bigr|+T^{\Gamma(t+\delta)}(p_i(t+\delta))~
\geq~\bigl|x-p_i(t+\delta)\bigr|+T^{\Gamma(t)}(p_i(t+\delta))\\[3mm]
\qquad\qquad \geq~|x-p_i(t)|+t~>~|x-q|+T^\Gamma(q).
\enda\eeq
According to (\ref{ch}), when the barrier increases from $\Gamma(t)$ to 
$\Gamma(t+\delta)$, the point $x$ is still  reached in minimum time by a trajectory
passing through $q$.
We conclude that, for $x\in V$, 
 $$\rho(t,x)~=~|x-q|\qquad\implies\qquad \rho(t+\delta,x)~=~|x-q|.$$
 In other words, the function $\rho(\cdot, x)$ cannot have a downward jump.
 However, it may well jump upward, from $|x-p_i(t)|$ to $|x-q|$.
\v
{\bf 4.}
Combining the previous steps {\bf 1-2-3}, for a.e.~time $t>0$ we obtain 
$${d\over dt}\int_{R_1} \rho(t,x)\, dx~\geq~- {\Hat T^2+\Hat T\over 2}\cdot \sum_i\,\bigl|\dot p_i(t)
\bigr|~=~
- {\Hat T^2+\Hat T\over 2}\cdot{d\over dt} m_1(\Gamma(t)).$$
  This proves (\ref{c1}).

\begin{figure}[htbp]
   \centering
 \includegraphics[width=0.6\textwidth]{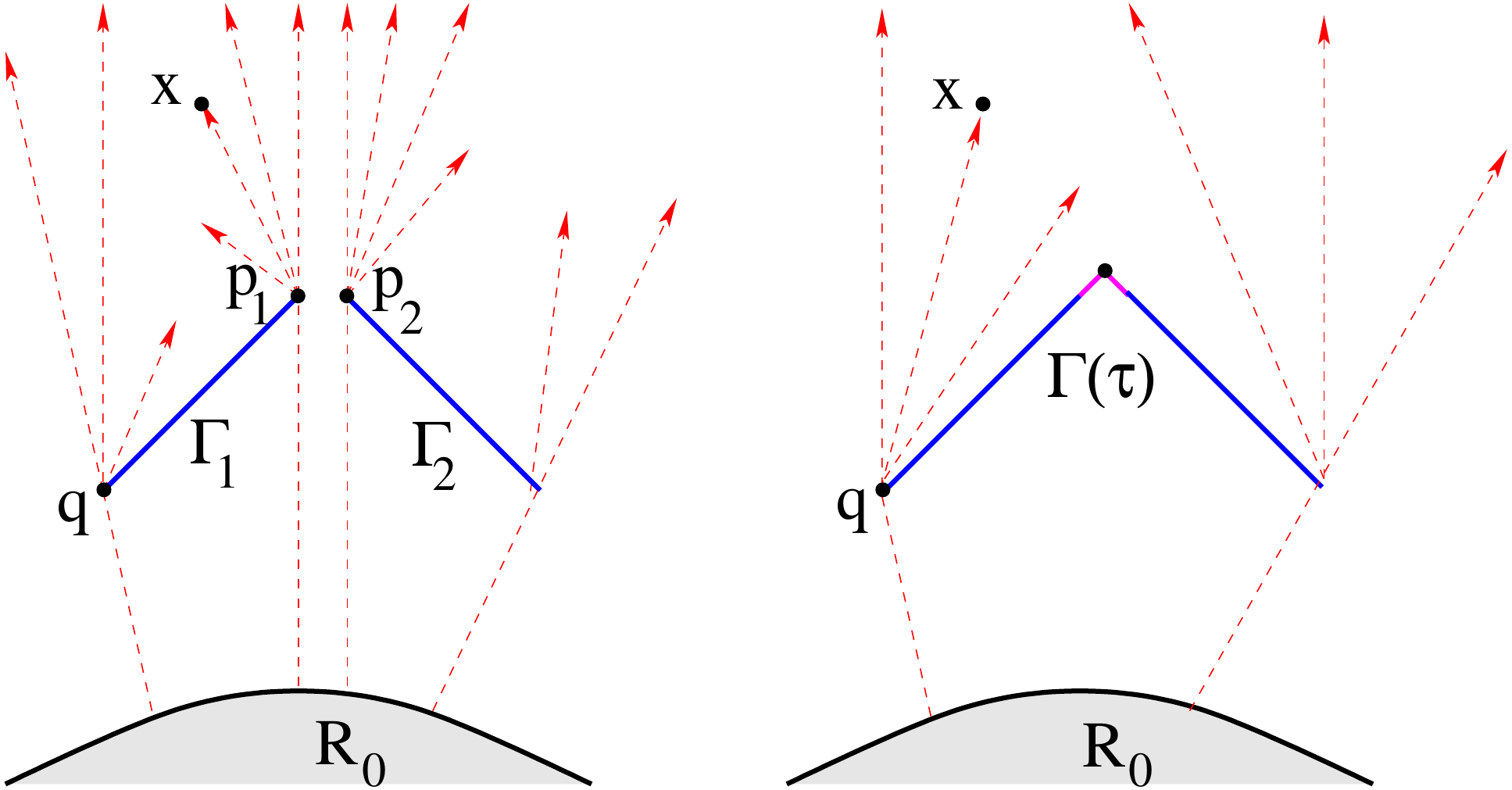}
   \caption{\small As time $t$ reaches a critical value $\tau$ when two 
   portions of the barrier join together, the topology of 
   $\Gamma(t)$ changes.  Both the minimum 
   time $T^{\Gamma(t)}(x)$ and the value $\rho(t,x)$ jump upward.}
   \label{f:fc148}
\end{figure}

To achieve the estimate (\ref{rb0}), we now observe that the integral in (\ref{c1}) depends continuously on time,
except at finitely many times $\tau_k$ where the topology
of $\Gamma(t)$ changes.  To understand what happens at these exceptional times,
as shown in Fig.~\ref{f:fc148}, left, assume that the  barrier $\Gamma(t)$ contains two segments $\Gamma_1$ and $\Gamma_2$ 
with moving endpoints $p_1(t)$, $p_2(t)$.
Assume that, at time $t=\tau$, the two segments join together: $p_1(\tau)=p_2(\tau)$
as in Fig.~\ref{f:fc148}, right.

Let $x\in R_1$ and assume that, for $t=\tau-\delta$ with $\delta>0$ small enough, 
the point $x$ is reached in minimum
time by a trajectory passing through $p_1(t)$. 
On the other hand, for $t=\tau$, assume that 
$\rho(\tau,x)~=~|x-q|$,
for some point $q$ along a different  optimal trajectory which reaches $x$ without crossing
$\Gamma(\tau)$.  
For $t<\tau$ we now have
\bel{rr1}\rho(t,x)=~|x-p_1(t)|~=~T^{\Gamma(t)}(x)- t,\eeq
while at time $\tau$
\bel{rr2}\rho(\tau,x)~=~|x-q|~=~T^{\Gamma(\tau)}(x) - T^\Gamma(q).\eeq
Observing that 
$$T^{\Gamma(\tau)}(x)~\geq~\lim_{t\to \tau-} T^{\Gamma(t)}(x),\qquad
\qquad T^\Gamma(q)~\leq~\tau,$$
by (\ref{rr1})-(\ref{rr2}) we conclude
\bel{rju}
\rho(\tau,x)~\geq~\lim_{t\to \tau-} \rho(t,x).\eeq
This shows that, at a time $\tau$  where the topology of the barrier $\Gamma(\cdot)$ changes, 
the function $\rho$ can only have upward jumps.

It remains to observe that, when $t=0$, one trivially has $\Gamma(0)=\emptyset$ and
$$
\int_{R_1} \rho(0,x)\, dx ~=~\int_{R_1}  d(x,R_0)\, dx.$$
Hence from (\ref{c1}) we  conclude (\ref{rb0}).

\v
{\bf 5.}
The previous analysis has established the estimate
(\ref{rb0}) in the case where the boundary of $R_0$ is a finite union
of circular arcs, and the barrier $\Gamma$ is the union of finitely many segments.
By an approximation argument, we shall extend the result to 
a general initial domain $R_0$ and a general barrier $\Gamma$.

As an intermediate step,  we show that the estimate (\ref{rb0}) holds for a general initial set $R_0$,
assuming that $\Gamma$ has finitely many
connected components: 
$\Gamma=\Gamma_1\cup\Gamma_2\cup\cdots\cup\Gamma_N$.    

Indeed, consider a sequence of open sets $(R_{0,n})_{n\geq 1}$ such that:
\begi
\item[(i)] The boundary 
of each $R_{0,n}$ is a finite union of circular arcs.
\item[(ii)]  As $n\to\infty$ the closures of these sets converge in the 
Hausdorff distance~\cite{AC}, namely
$d_H(\ov{R_{0,n}},\, \ov{R_0})\to 0$\endi
Moreover, for each $k=1,\ldots,N$, let $(\Gamma_{k,n})_{n\geq 1}$
be a sequence of compact connected sets such that:
\begi
\item[(iii)] Each $\Gamma_{k,n}$ is the union of finitely many segments.
\item[(iv)] $m_1(\Gamma_{k,n})\leq m_1(\Gamma_k)$.
\item[(v)] As $n\to\infty$ we have the convergence in the Hausdorff distance:
$d_H(\Gamma_{k,n},\Gamma_k)\to 0$.
\endi
For $x\in R_1$, let $\gamma_{x,n}$
be a polygonal line reaching $x$ in minimum time.
More precisely,  $\gamma_{x,n}$ minimizes $ m_1(\gamma_{x,n})$
among all polygonal lines connecting $x$ to some point $y\in \partial R_0$
without crossing the barrier $\Gamma_n=\cup_{k=1}^N \Gamma_{k,n}$.

We now observe that, for a.e.~point $x\in R_1$, the function $\rho$ defined at (\ref{rhodef})
satisfies
\bel{rusc}
\rho(x)~\geq~\limsup_{n\to \infty} \rho_n(x).\eeq
Indeed, we can parameterize every curve $\gamma_{x,n}$ by arc-length, say
$s\mapsto \gamma_{x,n}(s)$, with 
$$\gamma_{x,n}(0)~=~x,\qquad\qquad \gamma_{x,n}(m_1(\gamma_{x,n}))~\in~\partial R_{0,n}.$$
By taking a subsequence, we can assume the uniform convergence
$\gamma_{x,n}\to \gamma_x$ on every subinterval $[0,\ell]$ with $\ell<m_1(\gamma_x)$.
If now the derivatives $\dot \gamma_{x,n}$ are constant over some 
initial interval $[0,\bar s]$, the same is true of the derivative  $\dot \gamma_x$
of the limit function $\gamma_x$. This proves (\ref{rusc}).

{}From the inequality (\ref{ltn}) in Lemma~\ref{l:26} it follows
$$\Hat T~\doteq~\sup_{x\in R_1} T^\Gamma(x)~\geq~\limsup_{n\to\infty} 
\Hat T^{(n)}~\doteq ~\limsup_{n\to\infty}~
\sup_{x\in R_1} T^{\Gamma^{(n)}}(x).$$

In turn, since all functions $\rho_n$ are uniformly bounded, we  have
\bel{ls2}\bega{rl}\ds
\int_{R_1}\rho(x)\, dx &\ds \geq~\limsup_{n\to\infty} \int_{R_1} \rho_n(x)\, dx~\geq~
\limsup_{n\to\infty}\int_{R_1} d(x,R_{0,n})\, dx - {\Hat T_n+\Hat T_n^2\over 2}\, m_1(\Gamma_n)
\\[4mm]
&\ds
\geq~\int_{R_1} d(x,R_0)\, dx - {\Hat T+\Hat T^2\over 2}\, m_1(\Gamma).\enda\eeq
\v
{\bf 6.} Finally, we consider the general case where $\Gamma=\cup_{k\geq 1}
\Gamma_k$ is the union of countably many compact, connected components.
We call $\rho_\nu(\cdot)$ the map in (\ref{rhodef}), replacing 
$\Gamma$  
with a finite union $\Gamma_\nu\doteq\cup_{k=1}^\nu \Gamma_k$.

Thanks to Lemma~\ref{l:27}, the same argument used to prove (\ref{rusc}) now yields
\bel{rusc2}\rho(x)~\geq~\limsup_{\nu\to\infty} \rho_\nu(x)\qquad\qquad 
\hbox{for a.e.~}~x\in R_1\,.\eeq
Moreover, since $\Gamma_\nu\subset\Gamma$ for every $\nu\geq 1$, we trivially have
$$\Hat T~\doteq~\sup_{x\in R_1} T^\Gamma(x)~\geq~\sup_{x\in R_1} T^{\Gamma_\nu}(x)~\doteq~\Hat T_\nu\,.$$

By the previous steps, we already know that the estimate (\ref{rb0}) holds
for every $\rho_\nu$.   
Taking the limit as $\nu\to \infty$ and using Lemma~\ref{l:27}
we conclude
$$\bega{l}\ds\int_{R_1}\rho(x)\, dx~\geq~\limsup_{\nu\to \infty}~ \int_{R_1}\rho_\nu(x)\, dx~\geq~\int_{R_1} d(x,R_0)\, dx - \lim_{\nu\to \infty} {\Hat T_\nu^2+\Hat T_\nu\over 2}\, m_1(\Gamma)\\[4mm]
\qquad\qquad\ds \geq~\int_{R_1} d(x,R_0)\, dx - \lim_{\nu\to \infty} {\Hat T^2+\Hat T\over 2}\, m_1(\Gamma).\enda$$
This completes the proof of Lemma~\ref{l:5}.
\endproof
\v
\section{Avoiding barriers more efficiently} 
\label{s:5}
\setcounter{equation}{0}
As before, we assume that  $\R^2\setminus \Gamma$ is connected. By the analysis in Lemma~\ref{l:22},
if $p,q\notin\Gamma$, then for any $\ve>0$ we can connect these two points with a path that 
does not cross $\Gamma$ and has length $\leq |p-q|+(1+\ve) m_1(\Gamma)$.
Indeed, one can start with the segment having $p,q$ as endpoints, and then insert
detours to avoid crossing each connected component of $\Gamma$.

In this section we prove a sharper result.  Namely, if the barrier is sufficiently sparse, we can 
connect the two points $p,q$ with a path that avoids $\Gamma$ and has length
just slightly larger than $|p-q|$.  We begin by studying the case where $\Gamma$ is the union of finitely many (possibly intersecting) closed segments, then generalize.

\begin{lemma}\label{l:41}
In the $t$-$x$ plane, consider
a barrier $\Gamma$ consisting of
finitely many (possibly intersecting)  segments, none of which is parallel to the $x$-axis.
Assume that, for every $t>0$, the total length of the portion of $\Gamma$ contained
in the strip $[0,t]\times\R$ satisfies
\bel{bb}\psi(t)~\doteq~m_1 \Big(\Gamma\cap ([0,t]\times\R)\Big)~\leq ~
\sqrt 2\, \ve t\,,\qquad\qquad t\in [0,T],\eeq
for some  $0<\ve<1$.
Then there exists a continuous map
$\xi:[0,T]\mapsto\R$ with Lipschitz constant $\ve$,  which satisfies $\xi(0)=0$ and 
whose graph does not cross $\Gamma$. 
\end{lemma}
\v
\begin{figure}[htbp]
   \centering
 \includegraphics[width=0.45\textwidth]{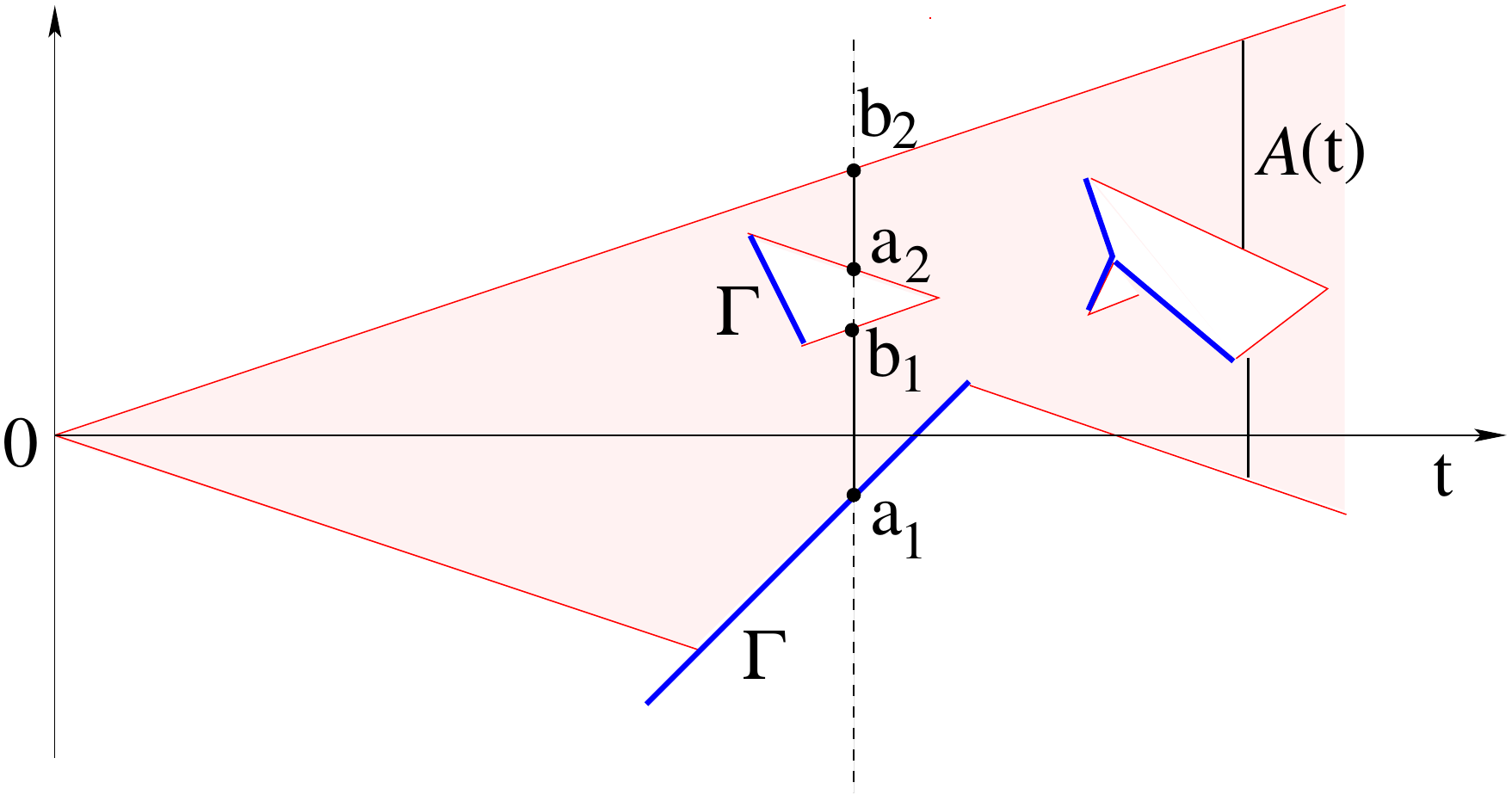}
 \quad \includegraphics[width=0.5\textwidth]{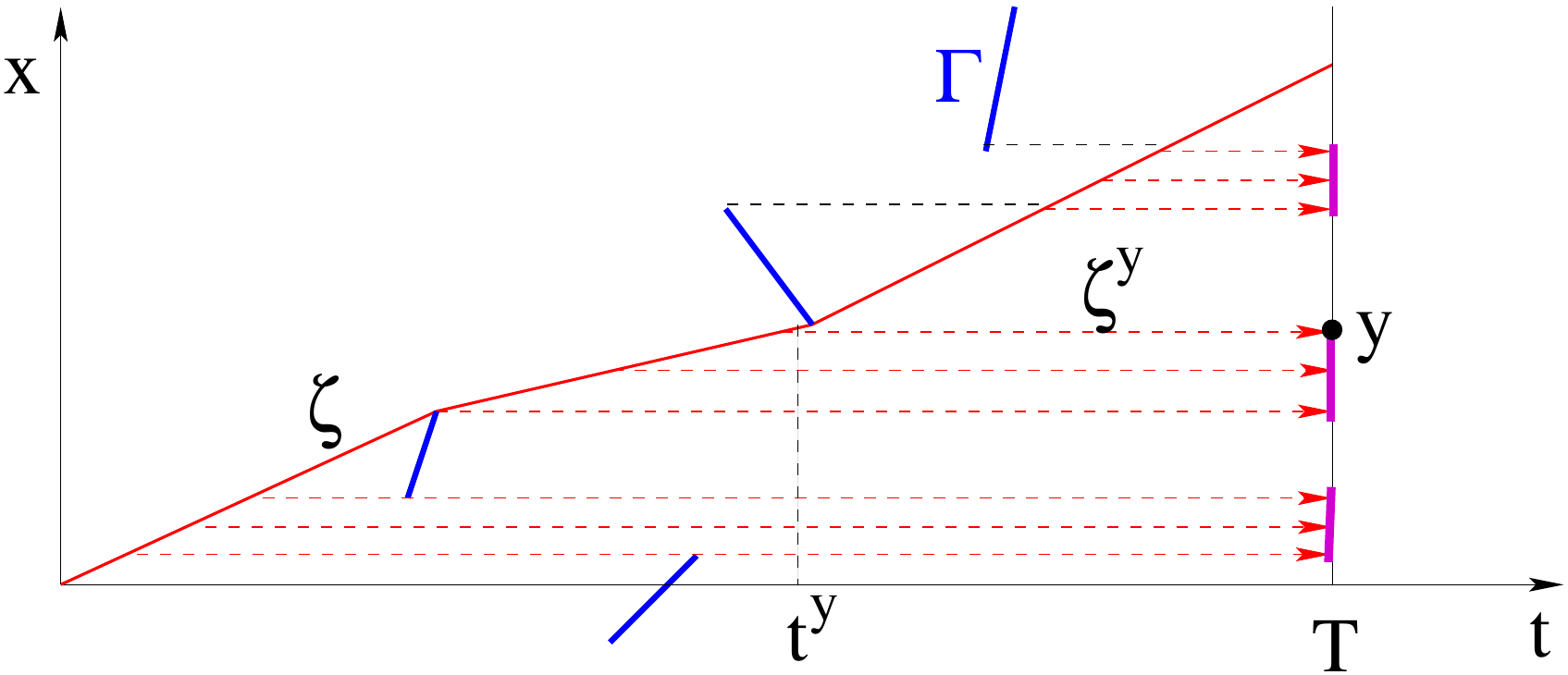}
   \caption{\small Left: for each $t>0$, the set $ \A(t)$ in (\ref{Rt})  is the union of finitely many segments $[a_k(t), b_k(t)]$. Right:
   the functions $\zeta$ and $\zeta^y$ constructed in the proof of Lemma~\ref{l:43}.
 }
   \label{f:fc178}
\end{figure}

{\bf Proof.} {\bf 1.} 
For every $t>0$, consider the set $\A(t)\subset\R$ of all values that can be attained by $\ve$-Lipschitz
functions, which are zero at the origin and whose graph does not cross $\Gamma$.
Namely, as shown in  Fig.~\ref{f:fc178}), left,
\bel{Rt}\bega{rl}
\A(t)&\doteq~\Big\{\xi(t)\,;~~\xi ~\hbox{is absolutely continuous},~~\xi(0)= 0,\\[3mm]
&\qquad \qquad \|\dot \xi\|_{\L^\infty}\leq \ve,~~
(s,\xi(s))\notin \Gamma ~~\hbox{for all} ~~s\in [0,t]\Big\}.\enda\eeq
Since $\Gamma$ is the union of finitely many segment,
we observe that each $\A(t)$ is the union of finitely many intervals, say 
$$\A(t)~=~\bigcup_k ~]a_k(t), b_k(t)[\,.$$
At any given time $t$, we denote by $\B(t)$ the set of the endpoints $a_k, b_k$ 
which lie along the barrier $\Gamma$, and by $\F(t)$ the set of the 
endpoints which are free, i.e. they do not lie on $\Gamma$.
The total length of the attainable set $\A(t)$ changes at the rate
\bel{rc}\bega{rl}\ds
{d\over dt} \Big(\meas(\A(t))\Big) &\ds=~\sum_k (\dot b_k(t) - \dot a_k(t))\\[4mm]
&\ds \geq~
\sum_{b_k(t)\in \F(t)} \dot b_k(t) - \sum_{a_k(t)\in\F(t)}  \dot a_k(t) 
-\sum_{b_k(t)\in \B(t)} |\dot b_k(t)| - \sum_{a_k(t)\in\B(t)}  |\dot a_k(t)|\,.
\enda
\eeq
On the other hand, 
from the definition of  $\psi$ at  (\ref{bb}), it follows
\bel{dpsi}\bega{rl}
\dot \psi(t)&\ds \geq~\sum_{a_k(t)\in\B(t)} \sqrt{ 1+ \dot a_k^2(t)} + \sum_{b_k(t)\in\B(t)} \sqrt{ 1+ \dot b_k^2(t)}\\[4mm]
&\ds \geq~\sum_{a_k(t)\in\B(t)}{1+|\dot a_k(t)|\over\sqrt 2} + \sum_{b_k(t)\in\B(t)} {1+|\dot b_k(t)|\over\sqrt 2}\,.\enda
\eeq
\v
{\bf 2.} 
To estimate the right hand side of (\ref{rc}), consider the function
\bel{fdef}f(t)~\doteq~\meas(\A(t)) -\sqrt 2(\ve t -\psi(t)).\eeq
By (\ref{bb}) it follows
$$\meas(\A(t))~=~f(t) +\sqrt 2 (\ve t -\psi(t)) ~\geq~f(t).$$
Therefore,
as long as $f(t)>0$, we have $\A(t)\not=\emptyset$.
In the remainder of the proof we will show that $f$ is positive and nondecreasing.

To begin, we observe that, for $t>0$ small, no barriers are present.   Hence
$$f(t) ~= ~\meas(\A(t))  -\sqrt 2\,\ve t~=~ 2\ve t-\sqrt 2 \,\ve t~ >~0.$$

Next, using (\ref{rc}) and (\ref{dpsi}), from (\ref{fdef}) we obtain
\bel{Aes}\bega{rl} \ds{d\over dt} f(t)&=~\ds {d\over dt} \meas(\A(t))- \sqrt 2\,\ve + \sqrt 2\cdot \dot 
\psi(t)\\[4mm]
&\ds\geq~\left(\#\F(t) \cdot \ve
-\sum_{b_k(t)\in \B(t)} |\dot b_k(t)| - \sum_{a_k(t)\in\B(t)}  |\dot a_k(t)|\right)
 - \sqrt 2\, \ve\\[4mm]
 &\qquad\ds
 + \sqrt 2\cdot 
  \left(\sum_{a_k(t)\in\B(t)}{1+|\dot a_k(t)|\over\sqrt 2} + \sum_{b_k(t)\in\B(t)} {1+|\dot b_k(t)|\over\sqrt 2}\right) \\[5mm]
  &\ds \geq~ \bigl(\#\F(t)+\#\B(t)\bigr)\cdot\ve   - \sqrt 2\, \ve ~\geq ~2\ve -\sqrt 2 \,\ve~>~0.\enda
\eeq
Here $\#\F$ and $\#\B$ denote the cardinality of the sets of free and constrained endpoints, respectively.
We observe that, as long as $\A(t)$ does not vanish,  its boundary contains at least two
points. This yields the last inequality in (\ref{Aes}).
\endproof
\v
Next, instead of (\ref{Rt}), we consider the sets
\bel{TA}\bega{rl}
\Tilde\A(t)&\doteq~\Big\{\xi(t)\,;~~\xi ~\hbox{is absolutely continuous},~~
 \dot \xi(s)\in [ \ve,\,3\ve]~~\hbox{for a.e.}~s\in [0,t]\,,\\[3mm]
&\qquad\qquad\qquad \qquad  \xi(0)= 0,\qquad 
(s,\xi(s))\notin \Gamma ~~\hbox{for all} ~~s\in [0,t]\Big\}.\enda\eeq
The same argument used to prove Lemma~\ref{l:41} yields

\begin{corollary}\label{c:42}
In the same setting as Lemma~\ref{l:41}, let (\ref{bb}) be replaced by 
\bel{bbb}\psi(t)~\doteq~m_1 \Big(\Gamma\cap ([0,t]\times\R)\Big)~\leq ~
{\sqrt 2\, \ve \over 1+2\ve}\,t\,,\qquad\qquad t\in [0,T].\eeq
Then, for every $t\in [0,T]$ the set $\Tilde A(t)$ in (\ref{TA}) is non-empty.
\end{corollary}

{\bf Proof.} Given a barrier $\Gamma\subset \R^2$ satisfying (\ref{bbb}), consider the shifted barrier
$$\Gamma^{2\ve}~=~\{ (t,x)\,;~~(t, x+2\ve t)\in \Gamma\}.$$
In view of (\ref{bbb}), this set satisfies the inequality 
\bel{bb3} \psi(t)~\doteq~m_1 \Big(\Gamma^{2\ve}\cap ([0,t]\times\R)\Big)~\leq ~(1+2\ve) 
m_1 \Big( \Gamma\cap ([0,t]\times\R)\Big)~\leq~
\sqrt 2\, \ve t\,,\qquad\qquad t\in [0,T].\eeq
Applying Lemma~\ref{l:41} we obtain an $\ve$-Lipschitz function 
$t\mapsto \xi(t)$ such that $\xi(0)=0$ and $(t,\Tilde \xi(t))\notin \Gamma^{2\ve}$
for all $t\in [0,T]$. 

Introducing the function $\Tilde\xi(t)~\doteq~\xi(t) + 2\ve t$,
we obtain $\Tilde \xi(t)\in \Tilde A(t)$ for all $t\in [0,T]$.   Hence $\Tilde \A(t)$ is nonempty.
\endproof
\v
In the next lemma, instead of (\ref{Rt}), for $t\in [0,T]$ we consider the attainable sets
\bel{A3} \bega{rl} \A_3(t)&\doteq~\Big\{\xi(t)\,;~~\xi ~\hbox{is absolutely continuous},~~
 \|\dot \xi(s)\|_{\L^\infty} \leq 3\ve\,,\\[3mm]
&\qquad\qquad\qquad \qquad  \xi(0)= 0,\qquad 
(s,\xi(s))\notin \Gamma ~~\hbox{for all} ~~s\in [0,t]\Big\}.\enda\eeq

\begin{lemma}\label{l:43}
In the $t$-$x$ plane, consider
a barrier $\Gamma$ consisting of
finitely many (possibly intersecting)  segments, none of which is parallel to the $x$-axis.
Assume that, for some  $0<\ve<1$,
\bel{bb4}\psi(t)~\doteq~m_1 \Big(\Gamma\cap ([0,t]\times\R)\Big)~\leq ~
{\ve\over 2}  t\,,\qquad\qquad \forall t\in [0,T].\eeq
Moreover assume that the total length of the barrier satisfies 
\bel{bb5}h~\doteq~m_1 (\Gamma)~\leq ~{\ve\, T\over 3}\,.\eeq
Then
\bel{A33}
m_1(\A_3(T)\cap [0,3h])~\geq~2h\,.\eeq
\end{lemma}

{\bf Proof.} {\bf 1.} Applying Corollary~\ref{c:42}, we obtain an absolutely continuous map $t\mapsto \zeta(t)$, 
with $\zeta(0)=0$, $\dot\zeta(t)\in [\ve, 3\ve]$, and whose graph does not intersect $\Gamma$.
\v
{\bf 2.}
Call  
$$V~\doteq~\bigl\{x\in\R\,;~~(t,x)\in \Gamma\quad\hbox{for some}~t\in [0,T]\bigr\}$$ 
the perpendicular projection of $\Gamma$ on the $x$-axis.
By (\ref{bb5}) it follows $m_1(V)\leq h$. Hence 
\bel{b8}m_1([0,3h]\setminus V)~\geq~2h\,.\eeq
\v{\bf 3.} Since $\dot \zeta(t)\geq\ve$, by (\ref{bb5}) 
for every $y\in [0,3h]\setminus V$, there exists a unique
time 
$t^y\in [0,T]$ such that $\zeta(t^y)=y$.
As shown in Fig.~\ref{f:fc178}, right, consider the map 
\bel{zyd}t~\mapsto~\zeta^y(t)~=~\left\{\bega{cl} \zeta(t)\qquad &\hbox{if}\quad t\in [0,t^y],\\[2mm]
y\qquad & \hbox{if}\quad t\in [t^y, T].\enda\right.\eeq
Our construction implies $(t,\zeta^y(t))\notin \Gamma$ for all $t\in [0,T]$.
Hence $\zeta^y(T)=y\in \A_3(T)$. 
We thus conclude
$$m_1(\A_3(T)\cap [0,3h])~\geq~
m_1([0,3h]\setminus V)~\geq~ 2h\,.$$
\endproof

\begin{remark}\label{r:44} {\rm Let $z=\zeta^y(\cdot)$ be one of the functions considered at (\ref{zyd}).
The length of its graph is computed by
\bel{lx}\bega{rl}
\ell &\ds=~\int_0^T\sqrt{ 1+ \dot z^2(t)}\, dt~\leq~\int_0^T \sqrt{ 1+ 3\ve \,\dot z(t)}\, dt~
\leq~\int_0^T\left(  1+ {3\ve\over 2}  \dot z(t)\right) dt \\[4mm]
&\ds
~\leq~T + {3\ve\over 2} z(T)~\leq~T + {3\ve\over 2}\cdot 3h ~=~T + {9\ve\over 2}\cdot m_1(\Gamma).
\enda
\eeq
This is a crucial bound, because it shows that the presence of a very sparse barrier
can lengthen the trajectories of the fire only  by an amount $\O(\ve)\cdot m_1(\Gamma)$.  
As a consequence, the time $\sigma^{-1}\cdot  m_1(\Gamma)$ spent for
constructing these walls is not compensated by
the additional time needed for the fire to go around them.}
\end{remark}

The final result proved in this section extends the previous lemmas to a general barrier 
$\Gamma=\cup_{i\geq 1}\Gamma_i$, which is the  union of countably many compact, connected,
rectifiable sets. As in Definition~\ref{d:adm}, we say that a  path
$t\mapsto \gamma(t)\in \R^2$, 
$t\in [0,\ell]$,  is {\bf admissible} if there exists a sequence of 1-Lipschitz paths
$t\mapsto \gamma_n(t)$ such that 
$\gamma_n(t)\notin \Gamma$ for all $t\geq 0$, and moreover
$\lim_{n\to\infty} \gamma_n(t)=\gamma(t)$,
uniformly for $t\in [0,\ell]$.
\v
\begin{lemma}\label{l:45} In the  $t$-$x$ plane, consider the points
$P=(-\kappa,0)$, $Q=(\kappa,0)$.
	Let $\Gamma\subset\R^2$ be a barrier such that, for every $r>0$,
	\bel{gp1}m_1\Big( \Gamma\cap \bigl([-\kappa, -\kappa+r]\times \R\Big)~<~{\ve\over 3}\,r,\qquad m_1\Big( \Gamma\cap \bigl([\kappa-r, \kappa]\times \R\Big)~<~{\ve\over 3}\,r.\eeq
	Moreover, assume 
	\bel{Gsm}h~\doteq~m_1(\Gamma)~<~ {2\,\kappa\,\ve\over 3 }\,.\eeq
	Then there exists an admissible path $\gamma: [0,\ell]\mapsto \R^2$
	such that 
	\bel{prop}\gamma(0)\,=\,P,\qquad \gamma(\ell)=Q,\eeq
	and with  length 
	\bel{glt} \ell~\leq~2\kappa + 
	9\ve \, m_1(\Gamma).\eeq
\end{lemma}
\v
{\bf Proof.}  {\bf 1.}
We begin by studying the case where $\Gamma= \cup_{i=1}^\nu\Gamma_i$ is the union of finitely many
compact, connected component. Then we will extend the result to the general case.

For any $\delta>0$, we can approximate each component $\Gamma_i$ with another connected set
$\Gamma_i'$, which is the union of finitely many closed segments, so that their Hausdorff distance
satisfies
\bel{HD} d_H(\Gamma_i,\Gamma'_i)~<~\delta\qquad\forall i=1,\ldots,\nu.\eeq
Moreover, we can assume that  (\ref{gp1}) still holds, with $\Gamma$ replaced by 
$\Gamma'= \cup_{i=1}^\nu\Gamma'_i$.

An application of Lemma~\ref{l:43}, with $[0,T]$ replaced by $[-\kappa, 0]$ yields 
the existence of a set $\A^-\subseteq [0, 3h]$ with the following properties.
\bel{pa1}m_1(\A^-)~\geq~2h,\eeq
For every $y\in \A^-$, there exists a Lipschitz function $t\mapsto \zeta^y(y)\in [0,y]$ such that
$$\dot \zeta^y(t)~\in~[0, 3\ve]\qquad  \hbox{for a.e.}~t\in [-\kappa, 0],$$
$$\zeta^y(-\kappa)=0,\qquad \zeta^y(0)= y,\qquad 
(t, \zeta^y(t))\notin\Gamma'\qquad\forall t\in [-\kappa,0].$$
Repeating the same argument on the interval $[0, \kappa]$, we obtain the existence of a set 
$\A^+\subseteq [0, 3h]$ such that 
\bel{pa2}m_1(\A^+)~\geq~2h,\eeq
For every $y\in \A^+$, there exists a Lipschitz function $t\mapsto \zeta^y(y)\in [0,y]$ such that
$$\dot \zeta^y(t)~\in~[-3\ve, 0]\qquad  \hbox{for a.e.}~t\in [0, \kappa],$$
$$\zeta^y(0)=y,\qquad \zeta^y(\kappa)= 0,\qquad 
(t, \zeta^y(t))\notin\Gamma'\qquad\forall t\in [0, \kappa].$$
By (\ref{pa1}) and (\ref{pa2}), we can choose $y\in \A^-\cap\A^+$.   Combining the two previous
constructions on $[-\kappa,0]$ and on $[0,\kappa]$, we obtain a Lipschitz function 
$\zeta:[-\kappa,\kappa]\mapsto [0, 3h]$ such that
$|\dot \zeta(t)|\leq 3\ve$ for a.e.~$t$, and  moreover
\bel{zpm}\zeta(0)\,=\,y,\qquad
\zeta(-\kappa)\,=\,\zeta(\kappa)\,=\,0,\qquad 
(t, \zeta(t))\notin\Gamma'\qquad\forall t\in [-\kappa,\kappa].\eeq
By the same argument used in Remark~\ref{r:44},  the length 
of the graph of $\zeta$ is bounded by
\bel{lgz}
\ell~=~
\int_{-\kappa}^\kappa \sqrt{ 1+ \dot \zeta^2(t)}\, dt~\leq  ~ 2\kappa+9\ve\cdot m_1(\Gamma).
\eeq

\begin{figure}[htbp]
   \centering
 \includegraphics[width=0.6\textwidth]{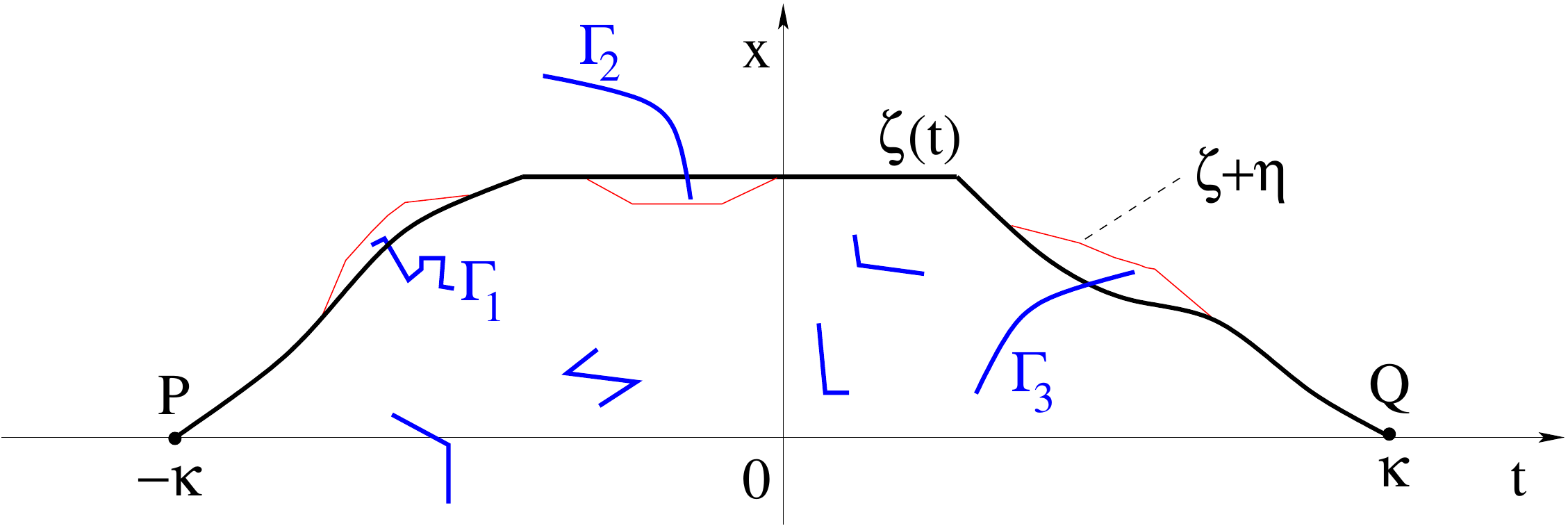}
   \caption{\small 
   Replacing the function $\zeta$ with $\zeta+\eta$, we obtain a new function 
   whose graph which does not intersect any of the connected components $\Gamma_1,\ldots,\Gamma_\nu$.
In this figure we have  $1,3\in I^-$ while $2\in I^+$.
 }
   \label{f:fc195}
\end{figure}
{\bf 2.} The graph of the function $\zeta$ constructed in the previous step does not touch
the components  $\Gamma'_i$, $1\leq i\leq\nu$ of the approximated barrier.
However,  it may well cross some components  $\Gamma_i$ of the original barrier.  
In this step (see Fig.~\ref{f:fc195}) we perform a small modification and construct a new map $z:[-\kappa, \kappa]\mapsto\R$
whose graph will not cross $\Gamma_1\cup\cdots\cup\Gamma_\nu$.

We begin by splitting 
$$\{1,\ldots,\nu\}~=~I^+\cup I^-,$$
where $I^+$ labels the components $\Gamma_i'$ lying above the graph of $\zeta$, 
while $I^-$ labels the components $\Gamma_i'$ lying below the graph of $\zeta$. 
We then set 
$$\Gamma~=~\Gamma^-\cup\Gamma^+,
\qquad \Gamma^-~\doteq ~\bigcup_{i\in I^-} \Gamma_i\,,\qquad  \Gamma^+~\doteq
~\bigcup_{i\in I^+} \Gamma_i
\,.$$

Consider the functions
$$\zeta^+(t) ~\doteq~\zeta(t)+2\delta,\qquad\qquad \zeta^-(t) ~\doteq~\zeta(t)-2\delta.$$
For any $\delta>0$ sufficiently small compared with $\ve$, by (\ref{HD}) the graph of $z^+$ does not intersect $\Gamma^-$,
while the graph of $z^-$ does not intersect $\Gamma^+$.
Call 
\bel{dnuu}\delta_\nu~\doteq~\min\Big\{ |p-q|\,;~~p\in \Gamma_i,~~q\in \Gamma_j\,,\quad 
1\leq i<j\leq\nu\Big\}~>~0,\eeq
and notice that
\bel{dnu} \min\Big\{ |p-q|\,;~~p\in \Gamma^-,~~q\in \Gamma^+\Big\}~\geq~\delta_\nu\,.\eeq

Consider the sets of times
$$\T^-~\doteq~\{ t\in [-\kappa, \kappa]\,;~~(t, x)\in \Gamma^+\quad\hbox{for some}~x<\zeta(t)\Big\},$$
$$\T^+~\doteq~\{ t\in [-\kappa, \kappa]\,;~~(t, x)\in \Gamma^-\quad\hbox{for some}~x>\zeta(t)\Big\}.$$
By (\ref{dnu}), there exists $\delta_\nu^*>0$ independent of $\delta$, such that 
\bel{dtt}\inf\Big\{ |t-t'|\,;~~t\in \T^+, ~~t'\in \T^-\Big\}~\geq ~\delta_\nu^*\,.\eeq
We now construct a Lipschitz function 
$$\eta:[-\kappa,\kappa]\mapsto [-2\delta, \, 2\delta]$$
such that 
$$\eta(-\kappa)\,=\,\eta(\kappa)~=~0,
\qquad \eta(t)~=~\left\{ \bega{rl} 2\delta\qquad &\hbox{if}~~t\in \T^+,\cr
-2\delta\qquad &\hbox{if}~~t\in \T^-.\enda\right.$$
By choosing $\delta>0$ sufficiently small, the Lipschitz constant of $\eta$ can be rendered as small as we like.   In particular, we can assume
$$\|\dot\eta\|_{\L^\infty}~\leq~2^{-\nu}.$$

The new function
$$z(t)\doteq \zeta(t) + \eta(t),\qquad t\in [-\kappa,\kappa],$$
has Lipschitz constant Lip$(z)\leq 3\ve+2^{-\nu}$.   Moreover, its graph does not intersect any
of the components $\Gamma_1,\ldots,\Gamma_\nu$.   Recalling (\ref{lgz}), 
the  length of the graph can be  bounded as
\bel{lze}\bega{rl}
\ell &\ds=~\int_{-\kappa}^\kappa\sqrt{ 1+ (\dot \zeta(t)+\dot \eta(t))^2}\, 
dt~\leq~\int_{-\kappa}^\kappa\sqrt{ 1+ \dot \zeta^2(t)}\, 
dt +\int_{-\kappa}^\kappa(\dot \zeta(t) + \dot \eta(t)) \dot \eta(t)\, 
dt\\[4mm]
&\ds \leq~2\kappa + 9\ve m_1(\Gamma) + 2\kappa (3\ve + 2^{-\nu}) 2^{-\nu}~\leq~2\kappa + 9\ve m_1(\Gamma) + 8\kappa\, 2^{-\nu}.
\enda
\eeq
\v
{\bf 3.} Next, consider the path $s\mapsto \gamma(s)$, $s\in [0,\ell]$, 
obtained by parameterizing the graph of 
$z$ by arc-length.    This is a 1-Lipschitz path that connects $P$ with $Q$, without touching any of the
connected components $\Gamma_1,\ldots, \Gamma_\nu$.   However, it may well 
cross many of the remaining components $\Gamma_i$, for $i>\nu$. 

To cope with this issue, we now use
Lemma~\ref{l:22} choosing $\epsilon= 2^{-\nu}$, and obtain a new path
$$\Tilde\gamma:[0, \Tilde\ell]\mapsto\R^2,$$
such that 
$$|\Tilde\gamma(0)-P|\,\leq \, 2^{-\nu},\qquad |\Tilde\gamma(\Tilde\ell)-Q|\,\leq \, 2^{-\nu},$$
$$\Tilde \gamma(s)\notin\Gamma\qquad\forall s\in [0,\Tilde\ell]\,.$$
The length of this new path is bounded by 
$$\Tilde\ell~\leq~\ell+  \sum_{i>\nu} m_1(\Gamma_i)
~\leq~ 2\kappa + 9\ve m_1(\Gamma) 
+ 8\kappa\, 2^{-\nu} +\sum_{i>\nu} m_1(\Gamma_i).$$
\v
{\bf 4.} 
By the previous steps, for every $\nu\geq 1$ there exists a  1-Lipschitz path
$$\gamma_\nu:[0,\ell_\nu]~\mapsto~\R^2\setminus \Gamma$$
such that
$$|\Tilde\gamma_\nu(0)-P|\,\leq \, 2^{-\nu},\qquad |\Tilde\gamma_\nu(\Tilde\ell)-Q|\,\leq \, 2^{-\nu}.$$
Moreover, its length satisfies
$$\ell_\nu~\leq~ 2\kappa + 9\ve m_1(\Gamma) 
+ 8\kappa \,2^{-\nu} + \sum_{i>\nu} m_1(\Gamma_i).$$
By Ascoli's compactness theorem, taking a subsequence we achieve the convergence
$\gamma_\nu\to \gamma$, where $\gamma:[0,\ell]\mapsto\R^2$ is a 1-Lipschitz path
joining $P$ with $Q$, with length
$$\ell~\leq~2\kappa + 9 \ve m_1(\Gamma).$$
By construction, this is an admissible path, satisfying the conclusion of the lemma.
\endproof

\begin{remark}\label{r:46}{\rm 
By the above construction, it follows that each path 
$\gamma_\nu$ differs by an amount $\O(1) \cdot \sum_{i>\nu}m_1(\Gamma_i)$ 
from the graph of a continuous function with Lipschitz constant $3\ve + 2^{-\nu}$.
Taking the limit, we thus obtain an admissible path $\gamma:[0,\ell]\mapsto \R^2$ which is the graph of a Lipschitz function $x= z(t)$, $t\in [-\kappa,\kappa]$
with Lipschitz constant $3\ve$.
}
\end{remark}

\begin{remark}\label{r:47} {\rm For simplicity, in the statements of Lemmas~\ref{l:41} and \ref{l:43}
we assumed 
a bound on the intersection of $\Gamma$ with the vertical strip $[0,T]\times \R$.   
Looking at the proofs, it is clear that 
we only needed a bound on the intersection of $\Gamma$ with the cone
$\{ (t,x)\,;~t\in [0,T], ~|x|\leq 3\ve t\}$.   In particular, the conclusion of Lemma~\ref{l:43}
remains valid if (\ref{bb4}) is replaced by 
\bel{bcone}
m_1\Big( \Gamma\cap \{ (t,x)\,;~t\in [0,T], ~|x|\leq 3\ve t\}\Big) ~\leq~{\ve\over 2} t\,,\qquad
\forall t\in [0,T].\eeq
The same remark applies to Lemma~\ref{l:45}.  Namely, all steps in the proof remain valid
if, for $ -\kappa<t<\kappa$, the assumption (\ref{gp1}) is replaced by
\bel{gp11}\bega{l}\ds
m_1\Big( \Gamma\cap \bigl\{ (t,x)\,;~|x|\leq 4\ve (t+\kappa)\bigr\}\Big)~<~{\ve\over 3}\,(t+\kappa),\\[3mm] \ds
m_1\Big( \Gamma\cap \bigl\{ (t,x)\,;~|x|\leq 4\ve (\kappa-t)\bigr\}\Big)~<~{\ve\over 3}\,(\kappa-t).
\enda
\eeq
}
\end{remark}

Thanks to the previous remarks, from Lemma~\ref{l:45} we deduce
\begin{corollary}\label{c:58}
Given $\theta_0,\ve_0>0$, there exists $\ve>0$ small enough so that the following holds.
Consider a triangle $\Delta_0$ with vertices
$$P=(-\kappa,0),\qquad Q=(\kappa,0),\qquad Z=(0, \theta_0\kappa).$$
Let $\Gamma\subset\R^2$ be a barrier such that, for every $r>0$,
\bel{bsp}
m_1\Big(\Gamma\cap\Delta_0\cap B(P,r)\Big)~\leq~\ve r,\qquad m_1\Big(\Gamma\cap\Delta_0\cap B(Q,r)\Big)~\leq~\ve r.\eeq
Then there exists a path $\xi:[0,\ell]\mapsto  \Delta_0$, joining $P$ with $Q$ without crossing the barrier $\Gamma$, with length bounded by
\bel{xle}
\ell~\leq~|P-Q| + \ve_0 \,m_1(\Gamma\cap\Delta_0).\eeq

\end{corollary}

\v
\section{Proof of Theorem~\ref{t:1}}
\label{s:6}
\setcounter{equation}{0}
Let $\Gamma$ be an optimal barrier for the optimization problem
{\bf (OP)}, and 
let $\Omega\subset\R^2$ be any open set.
We need to prove that the closure $\ov\Gamma$ does not contain all of $\Omega$.

Without loss of generality, we can assume $\Omega\subset\ov{R^\Gamma_\infty}$.
Otherwise, we can remove all barriers contained in the set
$\Omega\setminus \ov{R^\Gamma_\infty}$,
i.e., all portions of the wall which are never touched by the fire, and get a 
strictly smaller barrier. This yields a blocking strategy with a strictly lower cost.

As shown in Fig.~\ref{f:fc181}, right, the proof will be achieved by constructing 
a quadrilateral domain  $\Delta\subset\Omega$ with the following properties:
\begi
\item[(i)] The lower boundary $\gamma_0$ is the portion of a level set 
$\{ x\in\R^2\,;~T^\Gamma(x)=t_0\}$, between the two points $A$ and $B$.
\item[(ii)] For a suitable $h>0$, 
the upper boundary $\gamma^*$ is the portion of the curve
\bel{DDD}
\Big\{ x\in\R^2\,;~~
d(x,\gamma_0)~=~h\Big\}\eeq
between the points $C$ and $D$.
\item[(iii)]
The two sides $AC$ and $BD$ are segments which do not cross $\Gamma$,
and are part of optimal trajectories for the fire.
Their  lengths satisfy
\bel{ABCD}
|C-A|~=~|D-B|~=~d(C, \gamma_0)~=~d(D,\gamma_0)~=~h\,.\eeq
\item[(iv)] The total amount of barriers contained in $\Delta$ is small. Namely, for some
$\ve>0$ suitably small, one has
\bel{basm} m_1(\Gamma\cap\Delta)~\leq~\ve\,,\eeq
where $\sigma $ is the construction speed. Moreover, 
for every $s>0$ one has
\bel{sparse1}\bega{l}
\ds m_1\Big(\bigl\{y\in  \Gamma\cap\Delta\,; ~~d(y,\gamma_0)<s\bigr\}\Big)
~\leq~6 \ve s\,,\\[4mm]
\ds m_1\Big(\bigl\{y\in  \Gamma\cap\Delta\,; ~~d(y,\gamma^*)<s\bigr\}\Big)~
\leq~ 12\ve s\,.\enda\eeq
\endi
The first part of the proof, based on Lemma~\ref{l:5}, works out a construction
of the ``flow box" $\Delta$.
In the second part of the proof, using  Lemma~\ref{l:45}, we show that 
the reduced barrier
\bel{TGa}\Gamma^\diams~=~\Gamma\setminus \Delta\eeq
is still admissible, and yields a strictly lower cost.  We split the argument in several steps.
\v
{\bf 1.} Let $\ve>0$ be given.
Since the minimum time function $T^\Gamma$ is in SBV, it is differentiable at a.e.
point $\bar x$.   Moreover, the limit
\bel{rare}
\lim_{r\to 0+} {m_1(B(\bar x,r)\cap\Gamma)\over r^2}~=~0\eeq
also holds at a.e.~point $\bar x\in \Omega$.
We thus choose such a point $\bar x$, and  
consider a system of coordinates
with orthonormal basis $\{\bfe_1, \bfe_2\}$, where $\bfe_2= \nabla T^\Gamma(\bar x)$.
Call $\bar t = T^\Gamma(\bar x)$.

We now perform an affine transformation  of time and space coordinates, so that 
$(\bar t, \bar x)$ becomes the new origin of coordinates:
\bel{coc}
\Lambda(t', s_1, s_2) ~=~(\bar t + rt', ~\bar x + s_1 r \bfe_1 + s_2  r \bfe_2).\eeq
Since we are only interested in the local behavior of optimal trajectories for the fire 
in a neighborhood of $\bar x$, we consider a new problem where the initial
open set burned by the fire is 
$$R_0~\doteq ~\hbox{int} \bigl(\ov{R^\Gamma(\bar t)}\bigr).$$
Working in the $(t', s_1, s_2)$ coordinates, 
after renaming the variables and choosing a rescaling factor $r>0$ 
sufficiently small, we are led to study the following situation.
\begi
\item The total length of all barriers contained in the square $Q_2= [-2,2]\times [-2,2]$ satisfies
\bel{sh1}
m_1(\Gamma\cap Q_2)~\leq~\ve.\eeq
\item The initial set $R_0$ satisfies
\bel{sh2}\Big\{ (x_1, x_2)\in Q_2\,;~x_2\,< \,-\ve |x_1|\Big\}~\subseteq~
R_0\cap Q_2~\subseteq~\Big\{ (x_1, x_2)\in Q_2\,;~x_2\,<\,\ve |x_1|\Big\}.\eeq
\endi
\v
{\bf 2.} Call 
\bel{g13}
\Gamma_{1/3}~\doteq ~\Gamma\cap Q_2\cap\ov{R^\Gamma(1/3)}\eeq
 the portion of the barrier contained in the square $Q_2$ and 
touched by the fire within time $t=1/3$.  
By (\ref{sh1}) we trivially have
$m_1(\Gamma_{1/3})\leq\ve$.  

Our next goal is to apply Lemma~\ref{l:5} in this particular situation.
As in (\ref{R1def}), let  $R_1$ be the neighborhood of radius 1 around $R_0$.
For $x\in R_1$, define $\rho(x)$ as in (\ref{rhodef}),
with $\Gamma$ replaced by $\Gamma_{1/3}$.
By (\ref{gg}) and (\ref{HTB}) it now follows
\bel{ir3}\bega{l}\ds
\int_{R_1} \rho(x)\, dx ~\geq~\int_{R_1} d(x,R_0)\, dx - {\Hat T^2+\Hat T\over 2}\cdot m_1(\Gamma_{1/3})\\[4mm]
\qquad\qquad \ds\geq~\int_{R_1} d(x,R_0)\, dx - {(1+\ve)^2 + (1+\ve)\over 2}
\, \ve~\geq~\int_{R_1} d(x, R_0)\, dx - 2\ve,\enda\eeq
provided that $\ve>0$ is small enough.

\begin{figure}[htbp]
   \centering
 \includegraphics[width=0.80\textwidth]{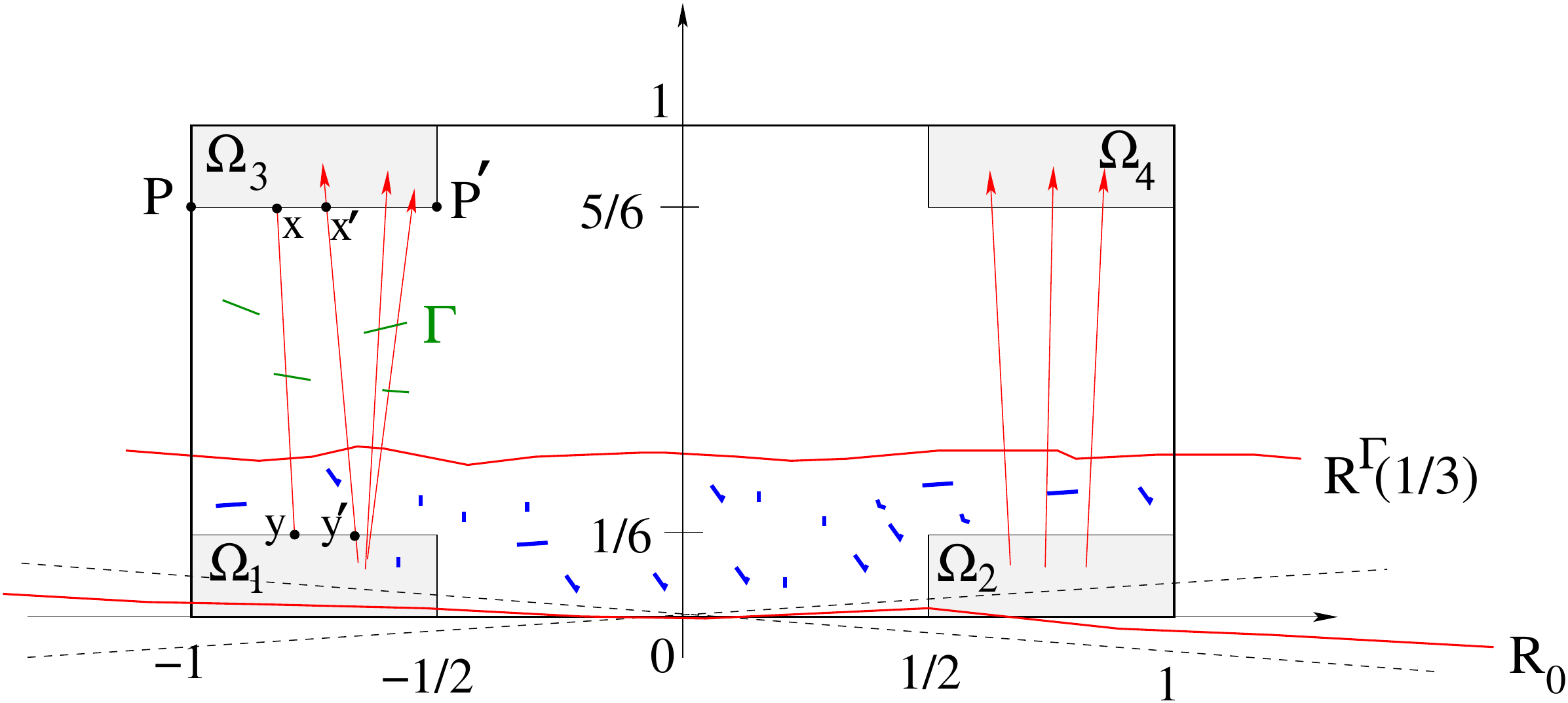}
   \caption{\small  Since the total length of all barriers is $\O(1)\cdot \ve$, 
   by Lemma~\ref{l:5} there exists
   many optimal trajectories that contain long straight segments, with one endpoint in $\Omega_1$ and the other in $\Omega_3$.   Since optimal trajectories do not cross each other, a barrier $\Gamma'\subset [-1,1]\times [1/4, 5/6]$,
   whose  total length is sufficiently small, cannot cross all of these segments.}
   \label{f:fc184}
\end{figure}

{}From (\ref{ir3}) we wish to conclude that,  within the square $Q_1=[-1,1]\times [-1,1]$,
most of the optimal trajectories for the fire contain long straight segments.
Since $\Gamma_{1/3}\subset B(0, 3)$, by the triangle inequality we have
\bel{id4}\rho(x)~=~d(x, R_0)\qquad\forall x\in R_1\setminus B(0,4).\eeq
By Lemma~\ref{l:22} it follows  
$$\rho(x)~\leq~T^{\Gamma_{1/3}}(x)~\leq~ d(x,R_0) + \ve,$$
and hence
\bel{int4}
\int_{R_1\cap (B(0,4)\setminus Q_1)}
 \rho(x)\, dx~\leq~
\int_{R_1\cap (B(0,4)\setminus Q_1)}  d(x, R_0) \, dx+ \ve\, m_2\Big(R_1\cap (B(0,4)\setminus Q_1)\Big).\eeq
{}From (\ref{ir3}), using (\ref{id4}) and then (\ref{int4}) we deduce
$$\left(\int_{R_1\cap  Q_1}+
\int_{R_1\cap (B(0,4)\setminus Q_1)}
\right) \rho(x)\, dx~\geq~\left(\int_{R_1\cap  Q_1}+
\int_{R_1\cap (B(0,4)\setminus Q_1)}
\right)  d(x, R_0) \, dx- 2\ve,$$
\bel{ir4}
\bega{rl}\ds\int_{R_1\cap  Q_1} \rho(x)\, dx&\ds\geq~\int_{R_1\cap  Q_1} d(x, R_0)\, dx
- \ve m_2\bigl(B(0,4)\setminus Q_1\bigr) - 2\ve\\[4mm]
&\ds\geq~\int_{R_1\cap  Q_1} d(x, R_0)\, dx - (16\pi+2) \ve.\enda\eeq
\v
{\bf 3.} As shown in Fig.~\ref{f:fc184},
consider in $Q_1$ the four rectangles 
$$\bega{l}\Omega_1~=~[-1, \,-1/2]\times [1,\, 1/6],\qquad 
\Omega_2~=~[1/2, \,1]\times [1,\, 1/6],\\[3mm]
\Omega_3~=~[-1, \,-1/2]\times [5/6, \,1],\qquad 
\Omega_4~=~[1/2, \,1]\times [5/6,\,1].\enda$$
Consider the lower side of $\Omega_3$. This is the horizontal segment $U$
with endpoints $P= (-1, 5/6)$ and $P' = (-1/2, 5/6)$.    By (\ref{ir4}), if $\ve>0$ 
is small enough, there exists a 1-dimensional subset $\Tilde U\subseteq U$ such that
$\rho(x)> 3/4$ for all $x\in \Tilde U$.   By choosing $\ve>0$ small, we can
make the size of $\tilde U$ as close to $1/2$ as we like.  Say,
$$m_1(\Tilde U)~>~1/3.$$

Given two distinct points  $x, x'\in \Tilde U$, let  $y, y'$ 
be the points where the optimal 
trajectories reaching $x,x'$ cross the boundary $\partial \ov{ R^{\Gamma}(1/3)}
= \{ z\in\R^2\,; ~T^\Gamma(z) = 1/3\}$.
Call $S,S'$ the segments with endpoints $x,y$ and $x', y'$, respectively.
Since these optimal trajectories are straight lines and   cannot cross each other
within their last segment  of length $3/4$,  we can find a constant $\lambda>0$ 
independent of $\ve$ such that 
\bel{baway} B\bigl(S, \lambda |x-x'|\bigr) \cap B\bigl(S', \lambda |x-x'|\bigr) ~=~\emptyset.\eeq
In other words, optimal trajectories reaching distinct points $x\in \Tilde U$ 
remain bounded away from each other.   A measure-theoretic argument now implies
that, if $m(\Gamma\cap Q_1)\leq \ve$ with $\ve>0$ small enough, we can find
at least one segment $S$ with endpoints $x,y$ as above, which does not intersect 
$\Gamma$.

\begin{figure}[htbp]
   \centering
 \includegraphics[width=0.8\textwidth]{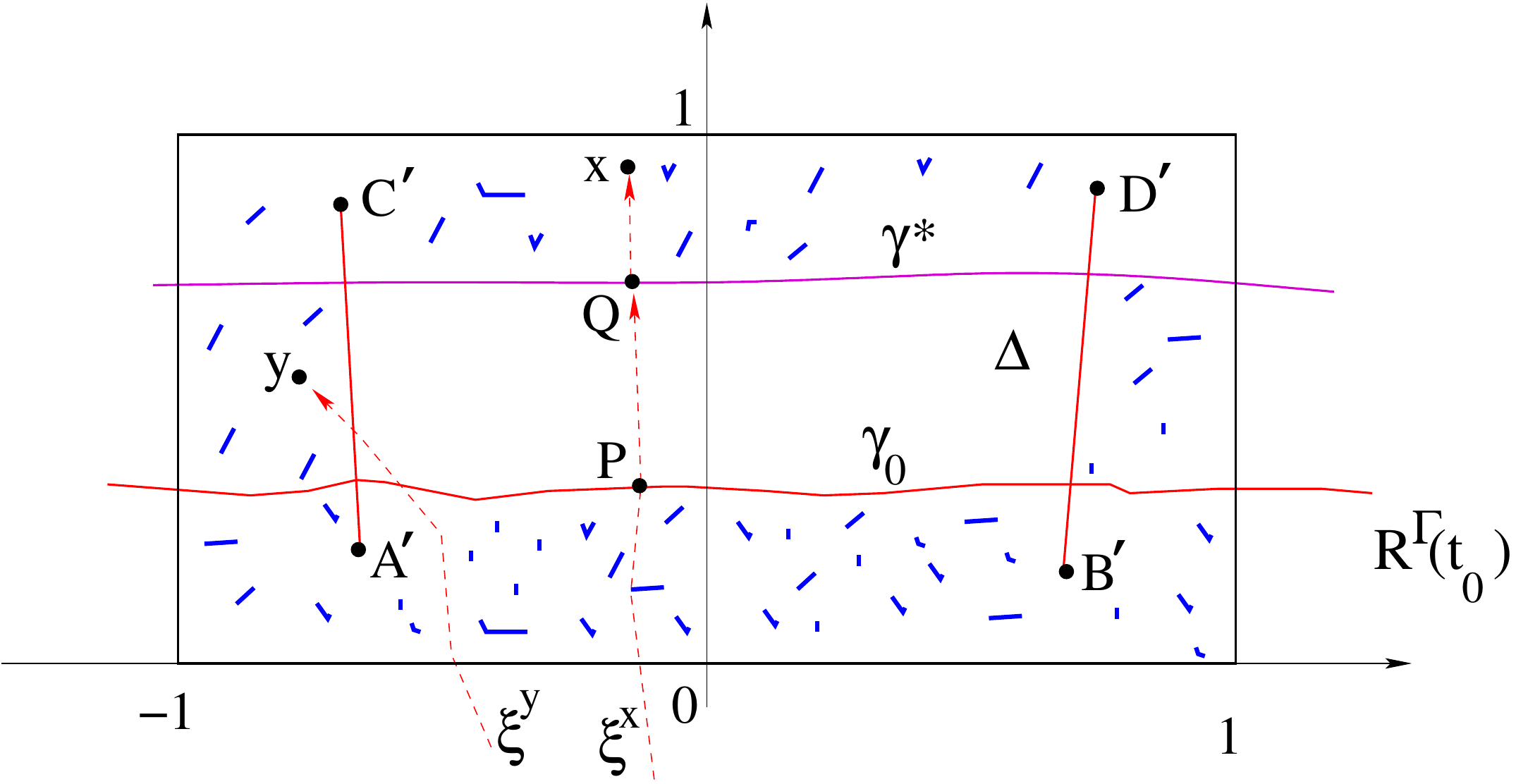}
   \caption{\small  If the domain $\Delta$ satisfies all properties (i)--(iv), optimal trajectories for the fire cannot exit from $\Delta$ through the 
   lateral boundaries $A' C'$ or $B' D'$.   In particular, $\xi^y$ is not optimal.
  }
   \label{f:fc199}
\end{figure}

\v
{\bf 4.} We are now ready to construct the quadrilateral domain $\Delta$ satisfying
the conditions (i)--(iv), as shown in Fig.~\ref{f:fc199}.

By the previous step, we can find four points 
$$A'\in \Omega_1\,,\quad B'\in \Omega_2\,,\quad C'\in \Omega_3\,,\quad D'\in \Omega_4\,,$$
with the following properties:
\begi
\item[(i)] When the barrier is taken to be $\Gamma_{1/3} \doteq
 \Gamma\cap \ov{R^\Gamma(1/3)}$,
the segment $A'C'$ is the last portion of a trajectory reaching $C'$ in minimum time.
Similarly, the segment $B'D'$ is the last portion of a trajectory reaching $D'$ in minimum time.
\item[(ii)] The segments $A'C'$ and $B'D'$ do not cross $\Gamma$.
\endi
As a consequence, for any $\tau\in [1/3, ~3/4]$, the segments $A'C'$ and $B'D'$ are still 
part of an 
optimal trajectory for the fire, in case the barrier $\Gamma$ is replaced by 
$\Gamma_\tau \doteq
 \Gamma\cap \ov{R^\Gamma(\tau)}$.

At this stage, it would be tempting to choose 
$\Delta$ as the quadrilateral having $\partial \ov{R^\Gamma(1/3)}$ as lower boundary,
the segments $A'C'$ and $B'D'$ as sides, and the curve
$$\gamma^*~=~\bigl\{x\in\R^2\,;~d(x ,\ov{R^\Gamma(1/3)}) = 1/3\bigr\}$$
as upper boundary. However, with this choice there is no guarantee that the 
bounds (\ref{sparse1}) will be satisfied.

To cope with this difficulty, 
the lower boundary will be chosen to be
$\gamma_0=\partial \ov{R^\Gamma(t_0)}$, for some $t_0\in [1/3, 1/2]$, while the upper boundary $\gamma^*$
will be the set of points having distance $h$ from the lower boundary, for some 
$h\in [1/4, 1/3]$.
The values of $t_0,h$ must be carefully chosen, in order to satisfy (\ref{sparse1}).

Consider the nondecreasing function 
$$\vp(t)~=~m_1\bigl(\Gamma\cap Q_2\cap\ov{R^\Gamma(t)}\setminus R_0\bigr).$$
By (\ref{sh1}),
$$\vp(0) ~\geq ~0,\qquad\qquad \vp(1)\leq \ve.$$
Using Riesz' sunrise lemma (see for example \cite{KF}, p.319)
we can find $t_0\in [1/3, 1/2]$
such that 
\bel{pic}
\vp(t_0+s)-\vp(t_0)~\leq~6\ve s ~\leq~{s\over 2}
\qquad\forall s\in [t_0,1]\,.\eeq
As in Lemma~\ref{l:29},
call $[a_i, b_i]$ the intervals during which the fire front touches the components $\Gamma_i$.
By (\ref{lsb}) we have
$$\bega{l}\ov{R^\Gamma(t_0+s)}
~\supseteq~B\Big(R^\Gamma(t_0),~ m_1\Big( [t_0, t_0+s]\setminus \bigcup_{i\geq 1}[a_i, b_i]\Big)\Big)\\[3mm]
~\qquad\qquad~\supseteq~B\bigl(R^\Gamma(t_0),~ s-\vp(t_0+s)+\vp(t_0)\bigr)~\supseteq~B\bigl(R^\Gamma(t_0), s/2\bigr).\enda
$$
In turn, by (\ref{pic}) this yields
$$\bega{rl}m_1\Big(\bigl\{y\in  \Gamma\cap\Delta\,; ~~d(y,\gamma_0)<s\bigr\}\Big)~
&\leq~m_1\Big( \Gamma\cap Q_1\cap \ov{R^\Gamma(t_0+2s)}~\setminus~
 \ov{R^\Gamma(t_0)}
\Big)\\[3mm]
&\leq~\vp(t_0+2s)-\vp(t_0)~\leq~12\ve\, s.\enda$$
Choosing $\ve>0$ small enough, this yields the first inequality in (\ref{sparse1}).

In a similar way, we now choose $h$ so that the second inequality in 
(\ref{sparse1}) is satisfied as well.   Consider the nondecreasing
function
$$\psi(t)~\doteq~m_1\Big( \Gamma\cap Q_1\cap ~\big\{ x\,;~d
\bigl(x, \ov{R^\Gamma(t_0)}\bigr)\leq t-t_0\bigr\}
\Big)$$
By (\ref{sh1}),
$$\psi(t_0) ~\geq ~0,\qquad\qquad \psi(1)\leq \ve.$$
Using Riesz' sunrise lemma
we can find $h\in [1/4, 1/3]$
such that 
\bel{picp}
\psi(t_0+h)-\psi(t)~\leq~12\ve (t_0+h-t) 
\qquad\forall t\in [t_0,\,t_0+h]\,.\eeq
Define the set
$$\gamma^*~\doteq~\bigl\{ x\in Q_1\,;~d(x, \ov{R^\Gamma(t_0)})~=~h\bigr\}.$$
By (\ref{picp}) it now follows
$$\bega{l}m_1\Big( \Gamma\cap Q_1\cap \bigl\{ x\,;~d(x,\gamma^*)\leq s\bigr\}\cap 
\bigl\{ x\,;~d\bigl(x, \ov{R^\Gamma(t_0)}\bigr)\leq h\bigr\}\Big)\\[3mm]
\qquad \leq~m_1\Big( \Gamma\cap Q_1\cap
\bigl\{ x\,;~d\bigl(x, \ov{R^\Gamma(t_0)}\bigr)\in [h-s, h]\bigr\}
\Big)\\[3mm]
\qquad\leq
~\psi(t_0+h) -\psi(t_0+h-s)~\leq~12\ve \,s.\enda $$
Choosing $\ve>0$ small enough, we thus obtain the second 
inequality in (\ref{sparse1}).

As shown in Fig.~\ref{f:fc199}, 
the quadrilateral domain $\Delta$ is now defined to be the set
of all points $x\in Q_1$ such that 
$$0~<~d\bigl(x,~ \ov{R^\Gamma(t_0)}\bigr)~<~h,$$
bounded between the two segments $A'C'$ and $B'D'$.

\v

{\bf 5.} Having constructed the domain $\Delta$, we now 
define the reduced barrier $\Gamma^\diams$ as in 
(\ref{TGa}), by removing all portions 
 inside $\Delta$.   Using the fact that $\Gamma$ is admissible,
we will show that $\Gamma^\diams$ is admissible as well. 
By (\ref{ad1}), this means
\bel{ad5} m_1\bigl(\Gamma^\diams\cap \ov{R^{\Gamma^\diams}(t)}\bigr)~\leq~\sigma t\qquad\forall t\geq 0.\eeq
For $t\leq t_0$ we trivially have
$$m_1\bigl(\Gamma^\diams\cap \ov{R^{\Gamma^\diams}(t)}\bigr)~=~
 m_1\bigl(\Gamma\cap \ov{R^{\Gamma}(t)}\bigr)~\leq~\sigma t .$$
For $t_0<t<  t_0+h$, we claim that
\bel{ad6} m_1\bigl(\Gamma^\diams\cap \ov{R^{\Gamma^\diams}(t)}\bigr)~\leq~
 m_1\bigl(\Gamma\cap \ov{R^{\Gamma}(t)}\bigr)~\leq\sigma t .\eeq
To prove (\ref{ad6}) we show that, for every $y\notin\Delta$, one has the implication
\bel{TTG}
T^\Gamma(y)~<~t_0+h\qquad\implies\qquad 
T^{\Gamma^\diams}(y)~=~T^{\Gamma}(y).\eeq
Indeed,  let $t\mapsto \xi^y (t)$
be an optimal trajectory for the fire, reaching $y$ in minimum time without crossing
the barrier $\Gamma^\diams$.
If $T^{\Gamma^\diams}(y)~<~T^{\Gamma}(y)$, then $\xi^y$ must cross some
barrier contained in $\Gamma\cap\Delta$.  As shown in Fig.~\ref{f:fc199},
this trajectory must be partly inside $\Delta$, then exit through one of the 
sides, either $A'C'$ or $B'D'$.  But this is impossible, because our construction
implies that both of these segments are part of optimal trajectories for the fire,
and two optimal trajectories cannot cross each other.  
For $t<t_0+h$, the inequality (\ref{ad6}) is an immediate consequence of 
(\ref{TTG}).

To achieve a bound valid for $t\geq t_0+h$, we claim that
\bel{st1}
\sup_{x\in\gamma^*} T^\Gamma(x)~\leq~t_0 +h + {1\over 2\sigma} m_1(\Gamma\cap \Delta).
\eeq
Indeed, consider any point $Q\in \gamma^*$, and let $P\in \gamma$ be a point such that
$d(P,Q)~=~d(P,\gamma)~ =~h$. 
Using  Corollary~\ref{c:58}, if $\ve>0$ was chosen sufficiently small, we can find a path $\xi:[0,\ell]\mapsto\Delta$,
joining $P$ with $Q$ without crossing the original barrier $\Gamma$, 
whose length satisfies
$$\ell ~\leq~ 
h + {1\over 2\sigma}\,m_1(\Gamma\cap \Delta).$$  This yields (\ref{st1}).
In turn,  for every $x\in \ov{R_\infty^\Gamma}$ with $T^\Gamma(x)\ge t_0+h$,
the inequality (\ref{st1}) implies
\bel{tdga}
T^{\Gamma}(x)~\leq~ T^{\Gamma^\diams}(x) + {1\over 2\sigma} \, m_1(\Gamma\cap\Delta).\eeq
Therefore
$$\ov {R^{\Gamma^\diams}(t)}~\subseteq~
\ov{ R^\Gamma\left(t+{1\over 2\sigma} m_1(\Gamma\cap \Delta)\right)}.$$
For any  $t\geq t_0+h$, the admissibility of $\Gamma$ now implies
\bel{extime}\bega{l}
\ds m_1\Big(\Gamma^\diams\cap \ov {R^{\Gamma^\diams}(t)}\Big)
~\leq~m_1\left( \Gamma\cap 
\ov{ R^\Gamma\left(t+{1\over 2\sigma} m_1(\Gamma\cap \Delta)\right)}\right)
-m_1(\Gamma\cap \Delta)\\[4mm]
\qquad\leq~\ds \sigma\cdot \left(t+{1\over 2\sigma}m_1(\Gamma\cap \Delta) \right)
-m_1(\Gamma\cap \Delta)~\leq~\sigma t-{1\over 2}\,m_1(\Gamma\cap \Delta) ,
\enda
\eeq
showing that the reduced barrier $\Gamma^\diams$ is admissible as well.
\v
{\bf 6.}
Since $\ov{R^{\Gamma^\diams}_\infty }= \ov{R^\Gamma_\infty}$, but $m_1(\Gamma^\diams)
<m_1(\Gamma)$, if $c_0>0$ we immediately conclude that the total cost of the strategy $\Gamma^\diams$
is strictly smaller:
$$m_2(R^{\Gamma^\diams}_\infty ) + c_0\,m_1(\Gamma^\diams)~<~m_2(R^\Gamma_\infty ) + c_0\,m_1(\Gamma).$$
This contradicts the optimality of $\Gamma$.

In the case $c_0=0$ we observe that, by (\ref{extime}), having removed  all barriers
contained inside $\Delta$, we are left with a little extra time: $( 2\sigma)^{-1} \, m_1(\Gamma\cap \Delta)$.
We can use this time to construct a circumference that forever shields a small disc
 from the fire.   More precisely, let $D_0$ be an open disc with radius $r_0$, so that
 the length of its boundary $\Gamma_0=\partial D_0$ satisfies 
 $$m_1(\Gamma_0)~=~2\pi r_0 ~\leq ~{1\over 2}\, m_1(\Gamma\cap \Delta).$$
We choose $D_0\subset \ov{R^\Gamma_\infty}$ so that 
$$\ov D_0\cap \ov{R^{\Gamma^\diams}(t_0+h)}~=~\emptyset.$$
In this way, the barrier
$\Gamma^*=\Gamma^\diams\cup \Gamma_0$ is still admissible.   The corresponding 
burned set satisfies
$$\ov{R^{\Gamma^*}_\infty}~\subseteq~\ov {R^\Gamma_\infty}\setminus D_0\,,$$ which has a strictly smaller area.
Again, this yields a contradiction with the optimality of $\Gamma$, proving the theorem.
\endproof

\end{document}